\documentclass[11pt]{article}

\usepackage[T1]{fontenc}
\usepackage[utf8]{inputenc}
\usepackage{lmodern}
\usepackage{microtype}
\usepackage{amsmath,amssymb,amsfonts}
\usepackage{mathtools}
\usepackage{graphicx}
\usepackage{booktabs}
\usepackage{siunitx}
\DeclareSIUnit{\AU}{AU}
\usepackage[round,authoryear]{natbib}
\usepackage{hyperref}
\usepackage{xcolor}
\usepackage[margin=1in]{geometry}

\title{ChromOps.jl: High-order simulation and discrete forward sensitivity analysis for chromatography models}

\author{Kristian Meyer$^{1}$\thanks{Corresponding author: Kristian Meyer
  (\href{mailto:KRVM@novonordisk.com}{KRVM@novonordisk.com}).}, Maksym Ratajczyk$^{1}$, Christopher Rackauckas$^{2,3}$\\[0.5em]
  $^{1}$Novo Nordisk A/S, Novo All\'e, Bagsv\ae rd, DK-2880, Denmark\\
  $^{2}$JuliaHub, Cambridge, MA 02138, USA\\
  $^{3}$Department of Mathematics, Massachusetts Institute of Technology, \\ Cambridge, MA 02139, USA}

\date{Preprint -- \today}

\begin{document}

\maketitle

\begin{abstract}
Mechanistic chromatography models are valuable tools in process systems engineering, but their use in gradient-based optimization can be hindered by computational cost and the effort required to obtain accurate objective-function gradients. To address this concern, a fully differentiable Julia chromatography solver, \texttt{ChromOps.jl}, is presented that combines high-order spatial discretization with discrete forward sensitivity analysis (DFSA). Two established high-order spatial discretizations, finite difference summation-by-parts (FD-SBP) and the discontinuous Galerkin spectral element method (DG-SEM), are compared on a 6-component ion-exchange chromatography problem with steric mass action kinetics. Both attain their theoretical convergence rates on manufactured problems, while FD-SBP offers better work-precision performance for both forward (primal) and dual-valued simulations. Moreover, its derivation and implementation require substantially less effort. The measured DFSA runtime follows $t_{\mathrm{DFSA}}\approx 0.9\,(1+N_p)\,t_{\mathrm{fwd}}$ for up to $N_p=24$ parameters propagated in a single \texttt{ForwardDiff.jl} chunk, corresponding to an incremental cost of approximately $0.9$ forward solves per additional parameter. By using the same code path for primal and dual-valued simulations, DFSA computes gradients of user-defined objectives through the discrete chromatography model without hand-derived sensitivity equations or chain-rule expressions. This reduces the implementation effort required to integrate mechanistic chromatography models into gradient-based optimization workflows in process systems engineering.
\end{abstract}


\section{Introduction}

Chromatography process development benefits from mechanistic modelling
tools \citep{Kumar2020MechanisticReview, Benner2019LabToManufacturing,
Rischawy2019GoodModelingPractice, Kobl2024OligoIEX}. Within process systems engineering (PSE), gradient-based optimization can be useful in e.g. parameter estimation and in-silico process optimization but two obstacles limit its routine use. First, each iteration requires both a forward solve and an accurate gradient of the objective with respect to many parameters. This combination is computationally expensive. Second, while current chromatography solvers \citep{Leweke2018CADET, Breuer2023DGSEMCADET, Meyer2020ChromaTech} provide parameter sensitivities, converting them -- via the chain rule -- into gradients of derived outputs (e.g., pool purity or yield) can be non-trivial and is left to the user. Addressing both obstacles requires a solver that is both fast in the forward pass and automatically differentiable for user-defined objectives. 

Mathematical models of chromatography processes
\citep{Guiochon2006Fundamentals, Nicoud2015ChromProcesses} consist of coupled
partial differential equations (PDEs) capturing transport, adsorption
and reaction phenomena within the chromatographic column. High-order
methods are attractive for solving such smooth (no discontinuities) PDEs because
their error decays rapidly under grid refinement, delivering high
accuracy with relatively few degrees-of-freedom (DOF). The discontinuous Galerkin spectral element method (DG-SEM) has been established as an efficient spatial discretization for chromatography \citep{Meyer2018NodalDGChromatography, Meyer2020ChromaTech, Breuer2023DGSEMCADET}. Its performance has been compared against second-order finite volume methods \citep{vonLieres2010CADET}, but comparisons against other high-order families are limited. 

An attractive alternative should retain DG-SEM’s high-order convergence and provable semi-discrete stability while being substantially cheaper to implement. Finite-difference summation-by-parts (FD-SBP) operators \citep{Svard2014SBPReview} satisfy these requirements. Like DG-SEM, they mimic integration by parts at the discrete level ensuring that their semi-discrete scheme is provably stable. The existence of diagonal-norm FD-SBP operators has been proved up to an interior order of $30$ \citep{AlbinKlarmann2016DiagonalNormSBP}. Their implementation is comparatively straightforward: unlike DG-SEM, it requires no variational formulation, basis functions, quadrature rules, element mappings, or numerical fluxes. Instead, it requires only an SBP differentiation operator and simultaneous-approximation-term (SAT) penalties for imposing boundary conditions. Two further structural differences distinguish the methods. FD-SBP does not duplicate interface nodes, as DG-SEM does, and can therefore represent a given spatial resolution with fewer DOF. Conversely, the upwind numerical flux used by DG-SEM introduces additional numerical dissipation that can be beneficial in under-resolved simulations \citep{Moura2015LinearDispersionDiffusion}. This effect is not exploited here since the presented chromatography models can be resolved directly at modest computational cost. The objective is to identify the simplest alternative to DG-SEM in terms of implementation effort. Therefore, other high-order method families, such as high-order finite-volume methods \citep{NORDSTROM2003453} or continuous spectral-element methods \citep{Xu2018SpectralHP}, are not considered, since their derivation and implementation are more involved than those of FD-SBP.

For sensitivity analysis, continuous forward sensitivity methods have
been developed for chromatography models \citep{Puttmann2013GRMSensitivities,
Meyer2023IndustrialIEXDG}, with \citet{Puttmann2016AlgorithmicDifferentiation}
employing automatic differentiation to assemble the Jacobians required
by the sensitivity equations. Adjoint methods have also been applied
to estimation and optimization \citep{Hahn2014AdjointChromatography}.
Discrete forward sensitivity analysis (DFSA), which differentiates the
discretized solver directly, has not been reported for chromatography models despite strong performance on broader differential-equation benchmarks relative to continuous methods
\citep{Ma2021ADSensitivityComparison}. Because DFSA automatically differentiates the discrete model together with post-processing, it yields sensitivities of arbitrary derived outputs directly, without requiring user-supplied chain-rule code. In principle, custom automatic differentiation rules can wrap continuous sensitivity methods to achieve the same composition, but existing chromatography implementations \citep{Leweke2018CADET, Breuer2023DGSEMCADET, Meyer2023IndustrialIEXDG} do not provide such rules. The computational cost of DFSA is evaluated by measuring the overhead relative to a primal solve, providing a direct measure of implementation efficiency. Combinations of DFSA with continuous sensitivity analysis for potential speedups are considered out of scope. 

Both obstacles -- computational cost and the need for user-supplied chain-rule code -- are addressed here in a single Julia \citep{Bezanson2017Julia}
framework, \texttt{ChromOps.jl}. Prior work \citep{Frandsen2025CADETJulia} used Julia as a rapid-prototyping platform for chromatography solvers; here it also serves as the production
language. Julia's just-in-time compilation via LLVM \citep{Lattner2004LLVM} lets numerical code run at speeds comparable to C or Fortran while remaining readable, so that a single code base can serve both research and production, eliminating the two-language problem.

The key contributions to chromatography PSE are:
\begin{itemize}
\item \texttt{ChromOps.jl}, a fully differentiable Julia chromatography
  solver developed as a Novo Nordisk inner-source project, providing a framework for gradient-based parameter estimation and in-silico process optimization with accurate gradients determined automatically;
\item a demonstration that DFSA yields gradients of arbitrary objective functions constrained by the discrete model at machine precision, at a per-parameter
  computational cost that scales linearly; and
\item a work-precision comparison of two established high-order
  discretizations, FD-SBP and DG-SEM, on a 6-component
  ion-exchange chromatography benchmark with steric mass action kinetics. This shows that
  FD-SBP can outperform DG-SEM under realistic preparative conditions while its derivation and implementation require less effort.
\end{itemize}

The paper is organized as follows. Section~\ref{sec:chromatography_model}
presents the lumped-rate chromatography model. Section~\ref{sec:methods}
covers the numerical methods: spatial discretization, time integration,
and DFSA. Section~\ref{sec:results} presents convergence studies on
manufactured problems, a 6-component ion-exchange benchmark comparing
spatial discretizations, linear solvers, ordinary differential equation (ODE) integrator performance,
and DFSA scaling. Section~\ref{sec:conclusion} summarizes the findings.

\section{The lumped-rate chromatography model}\label{sec:chromatography_model}
The chromatography column is modelled using a lumped-rate model \citep{Morbidelli1982POR} 
with axial dispersion and linear driving force mass transfer kinetics.
For each component $i \in \{s, 1, \ldots, N_c\}$, where $s$ denotes salt and $1, \ldots, N_c$ are the protein components, the governing equations are
\begin{align}
\frac{\partial c_i}{\partial t} &= -v_{\mathrm{int}}\frac{\partial c_i}{\partial z}
+ D_{\mathrm{ax},i}\frac{\partial^2 c_i}{\partial z^2}
- \frac{(1-\varepsilon)\varepsilon_p}{\varepsilon}k_{\mathrm{MT},i}(c_i-c_{p,i}), \\
\frac{\partial c_{p,i}}{\partial t} &= k_{\mathrm{MT},i}(c_i-c_{p,i})
- \frac{\partial q_i}{\partial t}.
\end{align}
Here, $c_i$, $c_{p,i}$, and $q_i$ denote bulk, pore, and adsorbed
concentrations, respectively, $v_{\mathrm{int}}$ is the
interstitial velocity, $D_{\mathrm{ax},i}$ is the axial dispersion coefficient,
$\varepsilon$ is the bulk (interstitial) porosity,
$\varepsilon_p$ is the particle porosity, and
$k_{\mathrm{MT},i}$ is the lumped mass-transfer coefficient between bulk and
pore phases. Moreover, $t\in[0,t_f]$ is the time domain with $t_f$ the final simulation time and $z\in\Omega=[0,L]$ is the spatial domain with $L$ being the chromatography column length.

At the column inlet and outlet, the boundary conditions of \citet{Danckwerts1953} are imposed,
\begin{subequations}\label{eq:danckwerts}
\begin{align}
v_{\mathrm{int}}c_{\mathrm{in},i}(t) &= v_{\mathrm{int}}c_i(t,0)
- D_{\mathrm{ax},i}\frac{\partial c_i}{\partial z}(t,0), \label{eq:danckwerts_in}\\
\frac{\partial c_i}{\partial z}(t,L) &= 0, \label{eq:danckwerts_out}
\end{align}
\end{subequations}
respectively. Here, $c_{\mathrm{in},i}(t)$ is the concentration entering the column. It results from the inlet concentration $c_{\mathrm{prog},i}(t)$ programmed into the chromatography system, typically defined by a piecewise linear polynomial, after passing through the dead volume and the mixer of the system. The dead volume $V_{\mathrm{dead}}$ is modelled as a perfect time delay $\tau_{\mathrm{dead}} = V_{\mathrm{dead}}/\dot{F}$, where $\dot{F}$ is the volumetric flow rate, and the mixer as a continuous stirred-tank reactor (CSTR) of volume $V_{\mathrm{mix}}$,
\begin{equation}\label{eq:cstr}
\frac{d c_{\mathrm{in},i}}{d t} = \frac{\dot{F}}{V_{\mathrm{mix}}}\left(c_{\mathrm{prog},i}(t-\tau_{\mathrm{dead}}) - c_{\mathrm{in},i}(t)\right),
\end{equation}
where $c_{\mathrm{prog},i}(t-\tau_{\mathrm{dead}}) = c_{\mathrm{prog},i}(0)$ for $t<\tau_{\mathrm{dead}}$, i.e.\ the dead volume is initially filled with the starting buffer.

The protein adsorption kinetics are described by the steric mass action (SMA) isotherm \citep{Brooks1992SMA}:
\begin{align}
\frac{\partial q_i}{\partial t} &= \bar{k}_{a,i}\left[c_{p,i}
- \frac{1}{k_{\mathrm{eq},i}}
\left(\frac{c_{p,s}}{\Lambda - \sum_{j=1}^{N_c}(\nu_j+\sigma_j)q_j}\right)^{\nu_i} q_i\right],
\qquad i = 1,\dots,N_c,
\end{align}
where $\bar{k}_{a,i}$ is the effective adsorption-rate coefficient, $k_{\mathrm{eq},i}$ is the equilibrium coefficient, $\Lambda$ is the ionic capacity, $\nu_i$ is the stoichiometric exchange coefficient, and $\sigma_i$ is the shielding coefficient. The corresponding salt balance follows as
\begin{equation}
\frac{\partial q_s}{\partial t} = -\sum_{i=1}^{N_c} \nu_i\frac{\partial q_i}{\partial t}.
\end{equation}
Since $q_s$ enters the model only through its time derivative, it is not carried as a state.

\section{Numerical methods}\label{sec:methods}

For both DG-SEM and FD-SBP, a method-of-lines strategy is followed: convection and
dispersion operators are derived individually, then combined into a single
semi-discrete spatial operator that is integrated in time
from appropriate initial conditions using an ODE integrator. Evaluating the solver using floating-point arithmetic yields the numerical solution, whereas evaluating it using dual-number arithmetic yields both the solution and its sensitivities with respect to the seeded parameters.

\subsection{Spatial discretization}

The continuous convection and dispersion operators are denoted
\begin{align}
\mathcal{L}_{\mathrm{conv}}(c_i) &= -v_{\mathrm{int}}\frac{\partial c_i}{\partial z}, \\
\mathcal{L}_{\mathrm{disp}}(c_i) &= D_{\mathrm{ax},i}\frac{\partial^2 c_i}{\partial z^2},
\end{align}
respectively. The boundary conditions are split between the two operators so
that their sum reproduces the Danckwerts conditions in Eq.~\eqref{eq:danckwerts}. The convection operator
imposes the inflow condition $c_i(t,0)=c_{\mathrm{in},i}(t)$ at $z=0$, while the dispersion operator imposes zero dispersive flux at both ends. At the outlet this gives the homogeneous
Neumann condition $\partial c_i/\partial z(t,L)=0$; at the inlet, the zero
dispersive flux combined with the convective feed recovers the Danckwerts
relation, leaving the inlet gradient determined by the flux balance rather than
fixed to zero. Each continuous operator is approximated by a discrete grid-based operator denoted with an $h$ superscript; the resulting semi-discrete system is
\begin{align}
\frac{d\mathbf{c}_i}{dt} &= \mathcal{L}^{h}_{\mathrm{conv}}(\mathbf{c}_i)
+ \mathcal{L}^{h}_{\mathrm{disp}}(\mathbf{c}_i)
- \frac{(1-\varepsilon)\varepsilon_p}{\varepsilon}
k_{\mathrm{MT},i}(\mathbf{c}_i-\mathbf{c}_{p,i}), \\
\frac{d\mathbf{c}_{p,i}}{dt} &= k_{\mathrm{MT},i}(\mathbf{c}_i-\mathbf{c}_{p,i})
- \frac{d\mathbf{q}_i}{dt}, \\
\frac{d\mathbf{q}_i}{dt} &= \bar{k}_{a,i}\left[\mathbf{c}_{p,i}
- \frac{1}{k_{\mathrm{eq},i}}
\left(\frac{\mathbf{c}_{p,s}}{\Lambda-\sum_{j=1}^{N_c}(\nu_j+\sigma_j)\mathbf{q}_j}\right)^{\nu_i}\mathbf{q}_i\right], 
\qquad i = 1,\dots,N_c, \\
\frac{d\mathbf{q}_s}{dt} &= -\sum_{i=1}^{N_c}\nu_i\frac{d\mathbf{q}_i}{dt},
\end{align}
where bold symbols denote vectors of nodal DOF.

For both FD-SBP and DG-SEM, the second-order dispersion operator
$\mathcal{L}^{h}_{\mathrm{disp}}$ is discretized using a Bassi-Rebay 1 (BR1) lifting
approach \citep{BassiRebay1997}. Rather than forming a single second-derivative
operator, an auxiliary gradient variable $\mathbf{g}_i=\partial_z\mathbf{c}_i$
is introduced and the first-derivative operator is applied twice: first to obtain $\mathbf{g}_i$, and then to obtain $\partial_z \mathbf{g}_i$.

\subsubsection{Finite difference summation-by-parts}

Given an ordered set of nodes $\{z_1,\dots,z_{N_i}\}$, where $N_i$ is the number of nodes, and an associated differentiation operator $D_{\textrm{FD}}$ with boundary order $r$ and interior order $l=2r$, the first derivative is approximated by
\begin{equation}
\dfrac{\partial \mathbf{c}_i}{ \partial z} \approx D_{\textrm{FD}}\mathbf{c}_i, 
\end{equation}
A positive definite norm (mass) matrix $H$ defines an $L^2$ inner product in the interval $[z_1, z_{N_i}]$, 
\begin{equation}
\int_{z_1}^{z_{N_i}} c_i v_j \, dz \approx \mathbf{c}_i^\top H \mathbf{v}_j,
\end{equation}
where $\mathbf{v}_j\in\mathbb{R}^{N_i}$ are values of a second function at the nodes. The pair $(D_{\textrm{FD}},H)$ satisfies the summation-by-parts property if 
\begin{equation}
H D_{\textrm{FD}} + D_{\textrm{FD}}^\top H^\top = B,
\end{equation}
where $B = \operatorname{diag}(-1, 0, \dots, 0, 1)$ is the boundary evaluation matrix. This identity is the discrete analogue of integration by parts and, together with suitable SAT penalties, yields an energy estimate that establishes stability of the resulting semi-discrete scheme in the discrete $L^2$ norm defined by $H$ \citep{Svard2014SBPReview}. The accuracy of the scheme is tuned by adjusting the number of nodes $N_i$, for a given boundary order $r$.

For convection, a Dirichlet SAT is imposed at the column inlet such that
\begin{equation}
\mathcal{L}^{h}_{\mathrm{conv}}(\mathbf{c}_i)
= -v_{\mathrm{int}}\left(D_{\textrm{FD}}\mathbf{c}_i
+ \boldsymbol{e}_1 H_{11}^{-1} \big(c_{i,1}-c_{\mathrm{in},i}\big)\right),
\end{equation}
where $\mathbf{e}_1=[1,0,\dots,0]^\top$ and $H_{11}$ is the first diagonal entry of the norm matrix.

For axial dispersion, the BR1 lifting approach described above is used, such that
\begin{align}
\mathbf{g}_i &= D_{\textrm{FD}}\mathbf{c}_i, \\
\mathcal{L}^{h}_{\mathrm{disp}}(\mathbf{c}_i) &= D_{\mathrm{ax},i}\left(
D_{\textrm{FD}}\mathbf{g}_i
+ \mathbf{e}_1 H_{11}^{-1} g_{i,1}
- \mathbf{e}_{N_i} H_{N_i N_i}^{-1} g_{i,N_i}\right),
\end{align}
where $g_{i,1}$ and $g_{i,N_i}$ are the boundary values of the auxiliary gradient and $\mathbf{e}_{N_i}=[0,\dots,0,1]^\top$. The last two terms are the SAT penalties that weakly set the numerical dispersive flux to zero at both column ends. At the outlet this is the homogeneous Neumann Danckwerts condition. At the inlet the vanishing dispersive flux and the
convective feed from $\mathcal{L}^{h}_{\mathrm{conv}}$ recover the Danckwerts
flux balance, setting the inlet gradient by that balance rather than to zero.

High-order FD-SBP operators are not unique, and multiple operators can be constructed for a given formal order \citep{Svard2014SBPReview}. Here, the minimum-stencil-width operators with a diagonal norm of \citet{MattssonAlmquist2018BoundaryOptimized} are used. The narrow stencils reduce the coupling between spatial DOF leading to a sparser state Jacobian, which can be advantageous for implicit time integration. These boundary-optimized operators introduce non-equispaced grid points near each boundary, while the interior grid remains equidistant. The node spacing near the boundaries is fixed by the operator, and these nodes provide additional free parameters that can be tuned to reduce truncation errors. The grid-point locations and operator coefficients are determined by solving a nonlinear constrained optimization problem. \citet{MattssonAlmquist2018BoundaryOptimized} derives the operators for boundary orders $r \in \{ 2,3,4,5,6 \}$, with corresponding interior orders $l=2r$. Here, the operators are constructed via \texttt{SummationByPartsOperators.jl} \citep{Ranocha2021SBPOperators} using its \texttt{LoopVectorization.jl} backend to enable SIMD vectorization. 

\subsubsection{Discontinuous Galerkin spectral element method}

The nodal Legendre--Gauss--Lobatto (LGL) DG-SEM formulation
\citep{HesthavenWarburton2008NodalDG} follows that of \texttt{Trixi.jl}
\citep{schlottkelakemper2021purely,ranocha2022adaptive}. The same implementation strategy of fusing multiply--add operations via \texttt{MuladdMacro.jl} is adopted to enable Fused Multiply-Add (FMA) instructions in the inner loops. The kernel is reimplemented here as a dedicated 1D convection--dispersion solver that is compatible with \texttt{ForwardDiff.jl} \citep{Revels2016ForwardDiff}.

The domain is partitioned into $N_e$ non-overlapping elements $\Omega_e$. On
each element, the solution is approximated by a polynomial of degree $p$,
\begin{equation}
c_{i}(t,z)\big|_{\Omega_e} \approx c_i^e (t, z) = \sum_{m=0}^{p} c_{i,m}^e(t)\,l_m(z),
\end{equation}
where $l_m$ is the Lagrange basis polynomial defined on the LGL nodes
$\{z_0,\dots,z_p\}$, satisfying $l_m(z_n) = \delta_{mn}$ with $\delta_{mn}$
the Kronecker delta. The coefficients $c_{i,m}^e(t)=c_{i}^e(t,z_m)$ are therefore the
nodal values of the solution on element $e$. The number of grid points is given by $N_i = N_e (p+1)$, i.e. there are two grid nodes at each element interface. The accuracy of the scheme is tuned by adjusting the number of elements $N_e$ for an arbitrary polynomial degree $p$.

The variational form for $\mathcal{L}^h_{\mathrm{conv}}$ is obtained by
multiplying the PDE by each Lagrange test function $l_n$ (the same basis used for the trial expansion of $c_{i}^e$), integrating over $\Omega_e$, and applying integration by parts. The resulting flux trace at $\partial\Omega_e$ is replaced by an upwind numerical
flux $v_{\mathrm{int}}\hat{c}_i$ to couple adjacent elements. This gives
\begin{equation}
\int_{\Omega_e} l_n\,\partial_t c_{i}^e\,dz
= v_{\mathrm{int}}\int_{\Omega_e}  l_n^{\prime} \,c_{i}^e\,dz
- v_{\mathrm{int}}\left[l_n\,\hat{c}_i\right]_{\partial \Omega_e},
\end{equation}
for $n=0,\dots,p$. In matrix form this becomes
\begin{equation}
M\frac{d\mathbf{c}_i^{(e)}}{dt}
= v_{\mathrm{int}} S^\top \mathbf{c}_i^{(e)}
- v_{\mathrm{int}}\left(\mathbf{e}_{p+1}^e \hat{c}_{i,R}
- \mathbf{e}_1^e \hat{c}_{i,L}\right),
\end{equation}
where $M_{nm}=\int_{\Omega_e} l_n l_m\,dz$ is the mass matrix, $S_{nm}=\int_{\Omega_e} l_n l_m ^{\prime} \,dz$ is the stiffness matrix, $\boldsymbol{c}_i^e \in \mathbb{R}^{p+1}$ is the element-local nodal state vector and $\mathbf{e}^e_1=[1,0,\dots,0]^\top$, $\mathbf{e}^e_{p+1}=[0,\dots,0,1]^\top$ are the first and last canonical basis vectors. The numerical fluxes $\hat{c}_{i,L}$ and $\hat{c}_{i,R}$ are the upwind traces at the left and right boundaries of element $e$, respectively. 

The integrals are approximated using LGL quadrature at the LGL
interpolation nodes, which is exact for polynomials of degree at most
$2p-1$. This collocation choice yields a diagonal mass matrix
$M_{nm}=\omega_n\delta_{nm}$, where $\omega_n$ are the LGL weights. This is similar to FD-SBP, whose norm matrix is also diagonal. The underlying integrand $l_m\,l_n$ reaches degree $2p$, so the mass matrix is inexact; the corresponding exact integration
gives a dense mass matrix \citep{HesthavenWarburton2008NodalDG}. However, as for the minimum-stencil-width FD-SBP operator, this makes the Jacobian more sparse, which can be beneficial for implicit time integration \citep{Breuer2023DGSEMCADET}. The stiffness matrix
integrand $l_n\,l_m^{\prime}$ has degree $2p-1$ and is therefore integrated
exactly by the same rule, giving $S^{\top}=D_{\textrm{DG}}^{\top}M$ with
$D_{\textrm{DG}, nm}=dl_m/dz(z_n)$ the nodal differentiation matrix. 

DG-SEM has a discrete summation-by-parts structure analogous to that of FD-SBP. On each element, the differentiation matrix $D_{\textrm{DG}}$ and diagonal LGL mass matrix $M$ satisfy
\begin{equation}
M D_{\textrm{DG}} + D_{\textrm{DG}}^\top M^\top = B,
\end{equation}
where $B$ is the same boundary-evaluation matrix as in FD-SBP. Thus, $M$ plays the same role in DG-SEM as the SBP norm matrix $H$ does in FD-SBP: each defines the discrete inner product used to formulate the summation-by-parts identity and derive the corresponding energy estimate.

Since the mass matrix $M$ is diagonal, inverting it is trivial:
\begin{equation}
\frac{d\mathbf{c}_i^{(e)}}{dt}
= v_{\mathrm{int}}\,M^{-1}S^\top\mathbf{c}_i^{(e)}
- v_{\mathrm{int}}\,M^{-1}\left(\mathbf{e}_{p+1}^e \hat{c}_{i,R} - \mathbf{e}_1^e \hat{c}_{i,L}\right).
\end{equation}

For the second-order operator $\mathcal{L}^h_{\mathrm{disp}}$, the BR1 lifting
described above is applied, repeating the same variational procedure twice:
\begin{align}
\mathbf{g}_i^{(e)} &= -M^{-1}S^\top\mathbf{c}_i^{(e)}
+ M^{-1}\left(\mathbf{e}_{p+1}^e \hat{c}_{i,R}
- \mathbf{e}_1^e \hat{c}_{i,L}\right), \\
\frac{d\mathbf{c}_i^{(e)}}{dt} &= - D_{\mathrm{ax},i} M^{-1}S^\top\mathbf{g}_i^{(e)}
+ D_{\mathrm{ax},i} M^{-1}\left(\mathbf{e}_{p+1}^e \hat{g}_{i,R}
- \mathbf{e}_1^e \hat{g}_{i,L}\right),
\end{align}
where $\hat{c}_i$ and $\hat{g}_i$ are central (interface-averaged) numerical
traces at interior element interfaces, with the boundary traces set by the
inflow and zero dispersive flux conditions described above.

\subsection{Time integration}

After spatial discretization, the coupled lumped-rate model is written as an
initial value problem (IVP),
\begin{equation}
\dot{\mathbf{y}}(t) = \mathbf{f}(\mathbf{y}(t), t; \boldsymbol{\theta}),
\qquad
\mathbf{y}(t_0) = \mathbf{y}_0,
\end{equation}
where $\mathbf{y} \in \mathbb{R}^N$ collects all the semi-discrete states (bulk, pore, and adsorbed concentrations), and $\boldsymbol{\theta} \in \mathbb{R}^{N_p}$ denotes model and operating parameters. This IVP is solved using the Julia ODE integrators provided by
\texttt{DifferentialEquations.jl} \citep{Rackauckas2017DiffEq}. Temporal errors are controlled through adaptive time stepping based on user-specified absolute and relative tolerances. 

The semi-discrete system is stiff, so implicit integration is used. Three
families are considered here: multistep backward differentiation formulas (BDF) \citep{Shampine1997MATLABODESuite},
diagonally implicit Runge--Kutta (DIRK) \citep{Kvaerno2004ESDIRK, KennedyCarpenter2019DIRKStiff}, and Rosenbrock--Wanner methods \citep{Steinebach2023Rodas5P}. BDF and DIRK methods solve the nonlinear equations arising from the implicit time discretization using Newton iterations,
\begin{equation}
W\,\Delta\mathbf{y}^{(m)} = -\mathbf{G}\!\left(\mathbf{y}^{(m)}\right),
\qquad
\mathbf{y}^{(m+1)} = \mathbf{y}^{(m)} + \Delta\mathbf{y}^{(m)},
\label{eq:newton_system}
\end{equation}
where $W$ is the system matrix, $m$ is the iteration index and $\mathbf{G}$ is the residual of the implicit time discretization. The system matrix is given by
\begin{equation}
W = I - \gamma h_n J,
\qquad
J = \frac{\partial\mathbf{f}}{\partial\mathbf{y}} \in \mathbb{R}^{N\times N},
\label{eq:iteration_matrix}
\end{equation}
where $\gamma$ is a coefficient set by the method, $h_n$ is the time step size at time step $n$, and $J$ is the state Jacobian. Rosenbrock--Wanner methods are linearly implicit and require just a single linear solve per stage, avoiding Newton iterations. This comes at the cost of having to update the system matrix $W$ at each time step; Rosenbrock variants exist that permit approximate Jacobians but are not considered here. The system matrix $W$ has the structural sparsity pattern of $J$, with the diagonal included. It is sparse, and non-symmetric, but its structural pattern remains fixed throughout integration. 

The linear system is solved by sparse LU factorization. The fixed structural sparsity pattern permits reuse of the symbolic analysis and fill-reducing ordering. Two operations therefore dominate the repeated linear-algebra work. Numerical refactorization recomputes the triangular factors $L$ and $U$ for an updated $W$, while each subsequent triangular solve applies these factors through forward and backward substitution. For BDF and DIRK methods, the factorization can be reused across multiple linear solves, with the frequency of refactorization determined by the solver. Rosenbrock--Wanner methods rebuild $W$ once per step, but all stages within that step share the same factorization.

\subsection{Discrete forward sensitivity analysis}

The goal is to compute the gradient of an objective function $\partial \phi(\mathbf{y})/\partial\theta_k$, where $\phi(\mathbf{y})$ is an arbitrary
objective derived from the solution (e.g., purity, yield). This gradient is required in numerical optimization when solved via a gradient-based algorithm.

The method relies on the first-order Taylor expansion of any differentiable function $f$:
\begin{equation}
f(a + b\varepsilon) = f(a) + f'(a)\,b\,\varepsilon,
\qquad \varepsilon^2 = 0.
\end{equation}
The nilpotent condition $\varepsilon^2=0$ truncates the series exactly after the linear term, so that derivative evaluation is exact (to machine precision) rather than approximate. A scalar carrying both the primal value and its derivative is called a dual number, $\hat{x}=x+\varepsilon x'$. Arithmetic on dual numbers is defined via operator overloading: each primitive operation ($+$, $\times$, $/$, $\sin$, $\exp$, etc.) has an overloaded method that applies the corresponding derivative rule, for example, $\hat{x}\cdot\hat{y} = xy + (x'y+xy')\varepsilon$, so the chain rule emerges automatically from the composition of these primitives.

\texttt{ForwardDiff.jl} \citep{Revels2016ForwardDiff} implements multi-dimensional dual numbers that propagate $N_\mathrm{chunk}$ partial derivatives simultaneously:
\begin{equation}
f\!\left(a + \sum_{k=1}^{N_\mathrm{chunk}} b_k\varepsilon_k\right)
= f(a) + f'(a)\sum_{k=1}^{N_\mathrm{chunk}} b_k\varepsilon_k,
\qquad \varepsilon_j\varepsilon_k = 0 \;\;\forall\, j,k.
\end{equation}
Here, $\varepsilon_1,\dots,\varepsilon_{N_\mathrm{chunk}}$ are independent nilpotent perturbation directions whose cross-products vanish, so that the partial derivatives do not interfere with one another. Each scalar in the program carries its $N_\mathrm{chunk}$ partials in a fixed-size tuple whose length is known at compile time. Because this length is part of the type, the compiler can unroll the per-partial arithmetic and emit efficient, vectorized machine code with no heap allocation. To compute sensitivities with respect to $N_\mathrm{chunk}$ parameters simultaneously, each parameter is \emph{seeded} with its own direction,
\begin{equation}
\hat{\theta}_k = \theta_k + \varepsilon_k, \qquad k = 1,\dots,N_\mathrm{chunk},
\end{equation}
and the chain rule propagates all $N_\mathrm{chunk}$ derivatives through every elementary operation in the semi-discrete model, time integrator, and post-processing. After propagation, extracting the $\varepsilon_k$ coefficient from any objective function $\phi$ yields the gradient $\partial \phi(\mathbf{y})/\partial\theta_k$. 

The chunk size is chosen to improve the memory access pattern and expose a fixed-width computation that the compiler can specialize for the selected chunk size, often yielding SIMD-vectorized arithmetic within each chunk; when the total number of parameters $N_p$ exceeds $N_\mathrm{chunk}$, the gradient is assembled over $\lceil N_p/N_\mathrm{chunk}\rceil$ passes, each propagating one chunk of $N_\mathrm{chunk}$ partials.

A practical consequence of differentiating the discrete solver is that no separate accuracy setting is needed for the sensitivities. Because the dual-valued solve traverses the same code path as the primal solve, the ODE integrator tolerances and the spatial resolution selected for the forward pass simultaneously control the accuracy of the computed sensitivities. However, on a given grid the sensitivities can be resolved less accurately than the primal solution, so the spatial resolution should be chosen with the sensitivities in mind.

\subsection{Software}

All benchmarks were run on a 64-bit x86-64 Windows 11 workstation with a 12th Gen Intel Core i7-1270P processor (Alder Lake; 4 performance cores and 8 efficiency cores; 2.20~GHz base clock) and 32~GB DDR5-5200 RAM. Experiments were executed under WSL2 using a single-threaded Julia v1.13.0 environment with a single-threaded BLAS. The process was pinned to a single performance core. The processor supports 256-bit AVX2 with FMA, but not AVX-512. For reproducibility, the main benchmark dependencies were: \begin{sloppypar}
\texttt{OrdinaryDiffEq.jl} v7.8.1, \texttt{OrdinaryDiffEqBDF.jl} v2.4.11, \texttt{OrdinaryDiffEqRosenbrock.jl} v2.7.4, \texttt{OrdinaryDiffEqSDIRK.jl} v2.9.6, \texttt{Sundials.jl} v6.7.1, \texttt{LinearSolve.jl} v5.18.0, \texttt{SummationByPartsOperators.jl} v0.5.101, and \texttt{ForwardDiff.jl}~\citep{Revels2016ForwardDiff} v1.4.6.
\end{sloppypar}

\section{Numerical examples}\label{sec:results}

This section verifies the high-order convergence of FD-SBP and DG-SEM on
manufactured problems with sufficient smoothness. The two spatial discretizations are benchmarked on a preparative 6-component ion-exchange problem. Finally, the cost scaling of DFSA is assessed.

All wall-clock times are reported as the minimum over at least seven passes, where each pass times every benchmark case once, so that repeated measurements of a case are interleaved with the other cases. Since timing measurements are typically right-skewed due to positive noise, the minimum reliably reflects the baseline runtime \citep{BenchmarkTools2016}.

\subsection{Sufficiently smooth manufactured problems}

Manufactured problems with sufficient smoothness and known analytical answers are used to verify that the implementation achieves its theoretical convergence rates, confirming correctness of the discrete operators and boundary treatment. Although practical chromatography problems are generally smooth in the sense that they contain no discontinuities, they may not possess the higher degree of smoothness required for the full theoretical convergence rates. Nevertheless, some of the benefits of the high-order discretization may still be retained.

For each spatial discretization and refinement level, the numerical
solution $c_h(t_f, \cdot)$ is compared against the exact solution at grid nodes. Reported
metrics are the $L^2$-error, the max-norm error ($L^\infty$), and the experimental order of
convergence (EOC),
\begin{equation}
\mathrm{EOC} = \frac{\log\left(e_1/e_2\right)}{\log(h_1/h_2)},
\end{equation}
where $e_1$ and $e_2$ are the errors on grids with spacing $h_1$ and $h_2$,
respectively. 

To keep temporal error from contaminating the spatial rates, time integration
uses tight tolerances ($\mathrm{abstol}=\mathrm{reltol}$), with \texttt{Vern9} at
$10^{-12}$ for the convection study and the implicit \texttt{QNDF} at $10^{-14}$
for the stiff dispersion study. Tightening these further left the errors
unchanged.

\subsubsection{Convection}

For convection, a sinusoidal convection test problem is considered
\begin{align}
\frac{\partial c}{\partial t} + v_{\mathrm{int}} \frac{\partial c}{\partial z} &= 0,
\qquad z\in[0,2\pi], \ t\in[0,t_f], \\
c(0, z) &= \sin(z), \\
c(t, 0) &= -\sin(v_{\mathrm{int}}t),
\end{align}
with constant speed $v_{\mathrm{int}}=2\pi$ and final time $t_f=2\pi$. The corresponding exact solution is
\begin{equation}
c(t, z) = \sin(z-v_{\mathrm{int}}t).
\end{equation}

The convergence results are reported in Table~\ref{tab:conv_convection}. Both DG-SEM and FD-SBP achieve their expected convergence rates of EOC $= p+1$ \citep{HesthavenWarburton2008NodalDG} and EOC $= r+1$ \citep{MattssonNordstrom2004, SvardNordstrom2006Accuracy}, respectively, confirming the correctness of the implementation. 

\begin{table}[t!]
\centering
\caption{History of convergence for the convection operator. The FD-SBP interior order is $l=2r$. FD-SBP is unavailable for $r=1$.}
\label{tab:conv_convection}
\begin{tabular}{@{}rrrrrr c rrrrrr@{}}
\toprule
\multicolumn{6}{c}{DG-SEM} & & \multicolumn{6}{c}{FD-SBP} \\
\cmidrule(lr){1-6}\cmidrule(lr){8-13}
$p$ & $N$ & \multicolumn{2}{c}{$L^2$} & \multicolumn{2}{c}{$L^\infty$} &
& $r$ & $N$ & \multicolumn{2}{c}{$L^2$} & \multicolumn{2}{c}{$L^\infty$} \\
\cmidrule(lr){3-4}\cmidrule(lr){5-6}\cmidrule(lr){10-11}\cmidrule(lr){12-13}
& & Error & EOC & Error & EOC & & & & Error & EOC & Error & EOC \\
\midrule
1 & 32  & 7.77e-02 & --   & 1.58e-01 & --   & & & & & & & \\
  & 64  & 2.00e-02 & 1.96 & 4.00e-02 & 1.98 & & & & & & & \\
  & 128 & 5.03e-03 & 1.99 & 1.00e-02 & 2.00 & & & & & & & \\
  & 256 & 1.26e-03 & 2.00 & 2.51e-03 & 2.00 & & & & & & & \\
\addlinespace
2 & 48  & 1.02e-03 & --   & 2.42e-03 & --   & & 2 & 256  & 6.77e-07 & --   & 1.08e-06 & --   \\
  & 96  & 1.30e-04 & 2.97 & 3.10e-04 & 2.97 & &   & 512  & 8.40e-08 & 3.01 & 1.32e-07 & 3.03 \\
  & 192 & 1.64e-05 & 2.98 & 3.91e-05 & 2.99 & &   & 1024 & 1.05e-08 & 3.01 & 1.63e-08 & 3.02 \\
  & 384 & 2.06e-06 & 2.99 & 4.91e-06 & 2.99 & &   & 2048 & 1.30e-09 & 3.00 & 2.02e-09 & 3.01 \\
\addlinespace
3 & 64  & 2.46e-05 & --   & 6.57e-05 & --   & & 3 & 64   & 1.07e-06 & --   & 1.95e-06 & --   \\
  & 128 & 1.54e-06 & 4.00 & 4.11e-06 & 4.00 & &   & 128  & 6.13e-08 & 4.13 & 1.12e-07 & 4.13 \\
  & 256 & 9.65e-08 & 4.00 & 2.58e-07 & 4.00 & &   & 256  & 3.68e-09 & 4.06 & 6.74e-09 & 4.06 \\
  & 512 & 6.03e-09 & 4.00 & 1.61e-08 & 4.00 & &   & 512  & 2.26e-10 & 4.03 & 4.14e-10 & 4.03 \\
\addlinespace
4 & 80  & 4.74e-07 & --   & 1.39e-06 & --   & & 4 & 32   & 2.89e-06 & --   & 6.04e-06 & --   \\
  & 160 & 1.48e-08 & 5.00 & 4.32e-08 & 5.00 & &   & 64   & 7.74e-08 & 5.22 & 1.71e-07 & 5.14 \\
  & 320 & 4.63e-10 & 5.00 & 1.36e-09 & 5.00 & &   & 128  & 2.24e-09 & 5.11 & 5.04e-09 & 5.08 \\
  & 640 & 1.45e-11 & 5.00 & 4.25e-11 & 4.99 & &   & 256  & 6.72e-11 & 5.06 & 1.52e-10 & 5.05 \\
\addlinespace
5 & 48  & 4.88e-07 & --   & 1.49e-06 & --   & & 5 & 32   & 2.53e-06 & --   & 7.57e-06 & --   \\
  & 96  & 7.67e-09 & 5.99 & 2.42e-08 & 5.94 & &   & 64   & 2.50e-08 & 6.67 & 9.58e-08 & 6.30 \\
  & 192 & 1.20e-10 & 6.00 & 3.77e-10 & 6.00 & &   & 128  & 2.46e-10 & 6.66 & 1.24e-09 & 6.27 \\
  & 384 & 1.88e-12 & 6.00 & 5.92e-12 & 5.99 & &   & 256  & 2.69e-12 & 6.52 & 1.70e-11 & 6.19 \\
\bottomrule
\end{tabular}
\end{table}

\begin{table}[t!]
\centering
\caption{History of convergence for the dispersion operator. The FD-SBP interior order is $l=2r$. FD-SBP is unavailable for $r=1$; for $r=5$, the minimum number of nodes is $N_i=17$.}
\label{tab:conv_dispersion}
\begin{tabular}{@{}rrrrrr c rrrrrr@{}}
\toprule
\multicolumn{6}{c}{DG-SEM} & & \multicolumn{6}{c}{FD-SBP} \\
\cmidrule(lr){1-6}\cmidrule(lr){8-13}
$p$ & $N$ & \multicolumn{2}{c}{$L^2$} & \multicolumn{2}{c}{$L^\infty$} &
& $r$ & $N$ & \multicolumn{2}{c}{$L^2$} & \multicolumn{2}{c}{$L^\infty$} \\
\cmidrule(lr){3-4}\cmidrule(lr){5-6}\cmidrule(lr){10-11}\cmidrule(lr){12-13}
& & Error & EOC & Error & EOC & & & & Error & EOC & Error & EOC \\
\midrule
1 & 16  & 2.61e-01 & --   & 2.86e-01 & --   & & & & & & & \\
  & 32  & 6.52e-02 & 2.00 & 7.09e-02 & 2.01 & & & & & & & \\
  & 64  & 1.63e-02 & 2.00 & 1.77e-02 & 2.00 & & & & & & & \\
  & 128 & 4.08e-03 & 2.00 & 4.42e-03 & 2.00 & & & & & & & \\
\addlinespace
2 & 24  & 6.30e-04 & --   & 9.29e-04 & --   & & 2 & 128  & 9.80e-08 & --   & 1.52e-07 & --   \\
  & 48  & 5.94e-05 & 3.41 & 9.81e-05 & 3.24 & &   & 256  & 6.01e-09 & 4.03 & 9.37e-09 & 4.02 \\
  & 96  & 6.83e-06 & 3.12 & 1.17e-05 & 3.07 & &   & 512  & 3.74e-10 & 4.01 & 5.83e-10 & 4.01 \\
  & 192 & 8.35e-07 & 3.03 & 1.44e-06 & 3.02 & &   & 1024 & 2.28e-11 & 4.03 & 3.62e-11 & 4.01 \\
\addlinespace
3 & 32  & 4.11e-05 & --   & 7.59e-05 & --   & & 3 & 64   & 7.35e-08 & --   & 1.47e-07 & --   \\
  & 64  & 5.06e-06 & 3.02 & 9.26e-06 & 3.03 & &   & 128  & 4.37e-09 & 4.07 & 8.40e-09 & 4.13 \\
  & 128 & 6.30e-07 & 3.01 & 1.15e-06 & 3.01 & &   & 256  & 2.65e-10 & 4.04 & 5.05e-10 & 4.06 \\
  & 256 & 7.87e-08 & 3.00 & 1.44e-07 & 3.00 & &   & 512  & 1.64e-11 & 4.02 & 3.19e-11 & 3.99 \\
\addlinespace
4 & 40  & 8.58e-07 & --   & 1.82e-06 & --   & & 4 & 16   & 4.79e-06 & --   & 9.38e-06 & --   \\
  & 80  & 2.43e-08 & 5.14 & 5.21e-08 & 5.13 & &   & 32   & 4.79e-08 & 6.65 & 9.05e-08 & 6.70 \\
  & 160 & 7.41e-10 & 5.03 & 1.59e-09 & 5.03 & &   & 64   & 6.04e-10 & 6.31 & 1.13e-09 & 6.32 \\
  & 320 & 2.31e-11 & 5.00 & 4.98e-11 & 5.00 & &   & 128  & 8.53e-12 & 6.15 & 1.62e-11 & 6.13 \\
\addlinespace
5 & 24  & 1.93e-06 & --   & 4.33e-06 & --   & & 5 & 16   & --       & --   & --       & --   \\
  & 48  & 5.68e-08 & 5.09 & 1.24e-07 & 5.13 & &   & 32   & 1.35e-07 & --   & 3.46e-07 & --   \\
  & 96  & 1.75e-09 & 5.02 & 3.78e-09 & 5.03 & &   & 64   & 1.09e-09 & 6.95 & 4.34e-09 & 6.32 \\
  & 192 & 5.46e-11 & 5.00 & 1.18e-10 & 5.00 & &   & 128  & 1.16e-11 & 6.55 & 6.37e-11 & 6.09 \\
\bottomrule
\end{tabular}
\end{table}

\subsubsection{Dispersion}

For the dispersion operator, the heat equation with homogeneous Neumann boundary conditions is used,
\begin{align}
\frac{\partial c}{\partial t} = D_{\textrm{ax}}\frac{\partial^2 c}{\partial z^2},&
\qquad z\in[0,L], \ t\in[0,t_f], \\
\frac{\partial c}{\partial z}(t, 0) = 0,& \quad
\frac{\partial c}{\partial z}(t, L) = 0,
\end{align}
with initial condition
\begin{equation}
c(0, z) = 100z(1-z).
\end{equation}
In the test setup, $L=1$, $D_{\textrm{ax}}=1/4$, and $t_f=0.5$. Using separation of variables, the exact solution can be written as
\begin{equation}
c(t, z)=\frac{a_0}{2}+\sum_{k=1}^{\infty}a_k
\cos\!\left(\frac{k\pi z}{L}\right)
\exp\!\left(-D_{\textrm{ax}}\frac{k^2\pi^2}{L^2}t\right),
\end{equation}
with
\begin{equation}
\frac{a_0}{2}=\frac{100}{6},
\qquad
a_k=-\frac{200\left((-1)^k+1\right)}{(k\pi/L)^2}.
\end{equation}
Only even modes are nonzero, so for $L=1$ the coefficients reduce to
$a_k=-400/(k\pi)^2$ for even $k$. In the numerical tests, the series is
truncated at $k=4$. Because odd modes are zero, the first omitted nonzero mode
is $k=6$, whose relative contribution is below $10^{-19}$. Thus, the truncated
series is already accurate to double precision.

Convergence results are reported in Table~\ref{tab:conv_dispersion}. DG-SEM matches the expected convergence rates of $p$ for odd polynomial degrees $\geq 3$, and $p+1$ for even degrees \citep{HesthavenWarburton2008NodalDG}. The case $p = 1$ reaches EOC of $p + 1$, one order above the periodic-domain lower bound of $p$ due to the non-periodic boundary treatment considered here. For FD-SBP, the observed EOC is $r+1$ for odd boundary orders and $r+2$ for even boundary orders, as expected \citep{MattssonNordstrom2004, SvardNordstrom2006Accuracy}.

\begin{table}[t!]
\centering
\caption{Model parameters for the 6-component IEX problem.}
\label{tab:chromatography_parameters}
\begin{tabular}{@{}llrl@{}}
\toprule
Parameter & Symbol & Value & Unit \\
\midrule
Length & $L$ & 10 & \unit{\milli\meter} \\
Mixer volume & $V_{\mathrm{mix}}$ & 0.1 & \unit{\milli\liter} \\
Dead volume & $V_{\mathrm{dead}}$ & 0.5 & \unit{\milli\liter} \\
Column volume & $V_{\mathrm{col}}$ & 1 & \unit{\milli\liter} \\
Bulk porosity & $\varepsilon$ & 0.37 & -- \\
Particle porosity & $\varepsilon_p$ & 0.66 & -- \\
Ionic capacity & $\Lambda$ & 0.3247 & \unit{\mol\per\liter} \\
Axial dispersion & $D_{\mathrm{ax}}$ & $1.5\times10^{-5}$ & \unit{\centi\meter\squared\per\second} \\
Mass transfer (salt) & $k_{\mathrm{MT,s}}$ & 1.39 & \unit{\per\second} \\
Mass transfer (proteins) & $k_{\mathrm{MT,i}}$ & 0.139 & \unit{\per\second} \\
Interstitial velocity & $v_{\mathrm{int}}$ & 0.0075 & \unit{\centi\meter\per\second} \\
\bottomrule
\end{tabular}
\end{table}

\begin{table}[t!]
\centering
\caption{SMA isotherm parameters for the 6-component IEX problem.}
\label{tab:isotherm_parameters}
\begin{tabular}{@{}llcccccc@{}}
\toprule
& & \multicolumn{6}{c}{Component} \\
\cmidrule(lr){3-8}
Parameter & Unit & A & B & C & D & E & F \\
\midrule
$\nu_i$ & -- & 22 & 22 & 21 & 20 & 10 & 23 \\
$\sigma_i$ & -- & 3 & 3 & 3 & 3 & 3 & 3 \\
$k_{\mathrm{eq},i}$ & -- & \num{e5} & \num{e5} & \num{3e4} & \num{5e2} & \num{5} & \num{e6} \\
$\bar{k}_{a,i}$ & \unit{\per\second} & 5.56 & 5.56 & 5.56 & 5.56 & 5.56 & 5.56 \\
\bottomrule
\end{tabular}
\end{table}

\subsection{Preparative ion-exchange chromatography with six components} 

In this benchmark, a 6-component ion-exchange (IEX) chromatography problem is considered using parameter values and operating conditions from an in-house Novo Nordisk A/S chromatography project. The model and isotherm parameters are given in Tables~\ref{tab:chromatography_parameters} and \ref{tab:isotherm_parameters}, respectively. The values are the production set point of the process, i.e.\ it represents a realistic manufacturing operating point.

To characterize the operating regime, the P\'eclet, Stanton, and Damk\"ohler numbers are defined using the volumetric flow rate per column volume $\dot{Q}=v_{\mathrm{int}}\varepsilon/L = 0.00278$~\unit{\per\second},
\begin{equation}
\mathrm{Pe} = \frac{v_{\mathrm{int}} L}{D_{\mathrm{ax}}} = 500,
\qquad
\mathrm{St}_i = \frac{k_{\mathrm{MT},i}}{\dot{Q}} =500\ (\text{salt}),\ 50\ (\text{proteins}),
\qquad
\mathrm{Da}_i = \frac{\bar{k}_{a,i}}{\dot{Q}} = 2000,
\label{eq:dimensionless_groups}
\end{equation}
which compare convection with axial dispersion, lumped mass transfer with convection, and adsorption kinetics with convection, respectively. 

\begin{figure}[t!]
\centering
\includegraphics[width=0.90\textwidth]{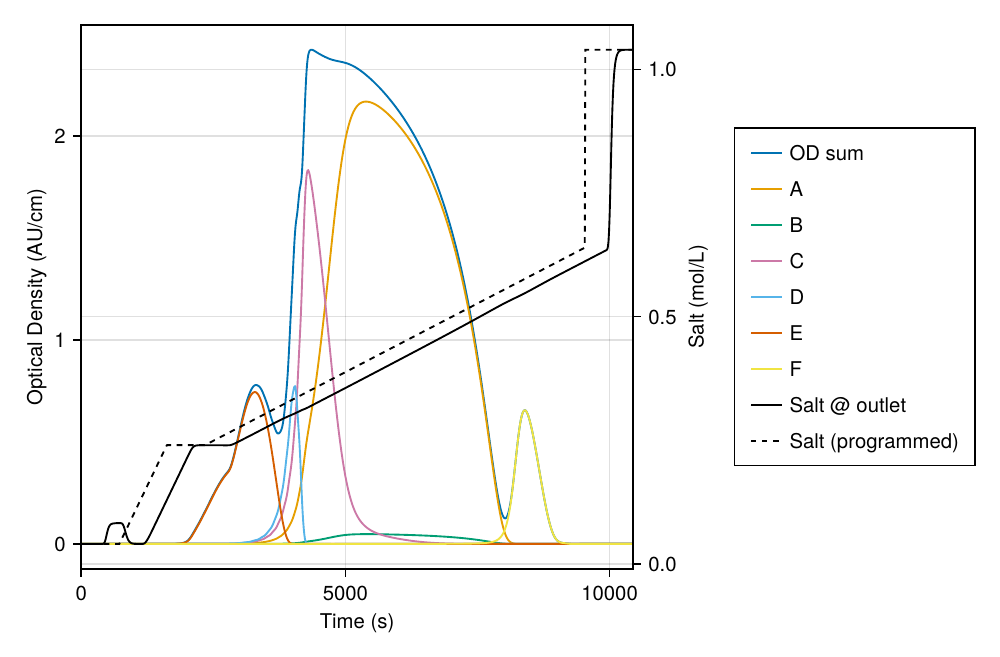}
\caption{Simulated chromatogram for the 6-component IEX problem.}
\label{fig:chromatogram}
\end{figure}

The large P\'eclet number $\mathrm{Pe} = 500$ indicates a convection-dominated regime. Although axial dispersion is small relative to convective transport, it remains finite and smooths the concentration fronts, keeping them steep but preventing shock-like discontinuities. The Stanton numbers -- $\mathrm{St}_i = 50$ for the proteins and $\mathrm{St}_s = 500$ for salt -- indicate that lumped mass transfer is fast relative to convection, although not instantaneous. Because the proteins are larger than salt, their mass transfer is slower; consequently, $\mathrm{St}_i$ is ten times smaller than $\mathrm{St}_s$. This slower protein mass transfer contributes to band broadening and motivates the use of a lumped-rate model. Moreover, the Damk\"ohler number $\mathrm{Da}_i = 2000$ indicates that adsorption is much faster than convection, with a negligible contribution to band broadening. Finally, the isotherm is strongly non-linear due to large characteristic charges $\nu_i\in[10,23]$. Hence, the problem is characterised by strong non-linearity, convection-dominated transport, and fast adsorption kinetics, typical of preparative ion-exchange chromatography.

The programmed inlet profile is a rectangular injection for the proteins:
\begin{equation}
c_{\mathrm{prog},i}(t) = \begin{cases}
c_{\mathrm{feed}}, & t \in [0,\, 360]~\unit{\second}, \\
0, & \text{otherwise},
\end{cases}
\end{equation}
where the feed concentrations are
\begin{equation}
\mathbf{c}_{\mathrm{feed}} = 10^{-5}\begin{bmatrix} 138.1 & 3.046 & 30.87 & 6.092 & 16.65 & 8.326 \end{bmatrix}^\top~\unit{\mol\per\liter}.
\end{equation}
For salt, the programmed inlet concentration follows a piecewise-linear gradient program:
\begin{equation}
c_{\mathrm{prog,s}}(t) = c_{\mathrm{s,start}}^{(k)} + \frac{c_{\mathrm{s,end}}^{(k)} - c_{\mathrm{s,start}}^{(k)}}{\Delta t^{(k)}}(t - t^{(k)}),
\qquad t \in [t^{(k)},\, t^{(k)}+\Delta t^{(k)}],
\end{equation}
where $k$ indexes the process phase and $c_{\mathrm{s,start}}^{(k)}$, $c_{\mathrm{s,end}}^{(k)}$, $\Delta t^{(k)}$ are the start concentration, end concentration, and duration of phase~$k$, respectively. The process consists of six phases with durations, start and end salt concentrations:
\begin{align}
\boldsymbol{\Delta t} &= \begin{bmatrix} 360 & 360 & 900 & 720 & 7200 & 720 \end{bmatrix}^\top~\unit{\second}, \\
\mathbf{c}_{\mathrm{s,start}} &= \begin{bmatrix} 0.04 & 0.04 & 0.04 & 0.24 & 0.24 & 1.04 \end{bmatrix}^\top~\unit{\mol\per\liter}, \\
\mathbf{c}_{\mathrm{s,end}} &= \begin{bmatrix} 0.04 & 0.04 & 0.24 & 0.24 & 0.64 & 1.04 \end{bmatrix}^\top~\unit{\mol\per\liter},
\end{align}
respectively. Outlet concentrations are converted to optical density (OD) with extinction coefficient $w = \qty{1.167e4}{\AU\liter\per\mol\per\centi\meter}$; \unit{\AU} denotes the dimensionless absorbance unit, and OD has dimensions of inverse length \unit{\AU\per\centi\meter}. 

\begin{figure}[t!]
\centering
\includegraphics[width=0.90\textwidth]{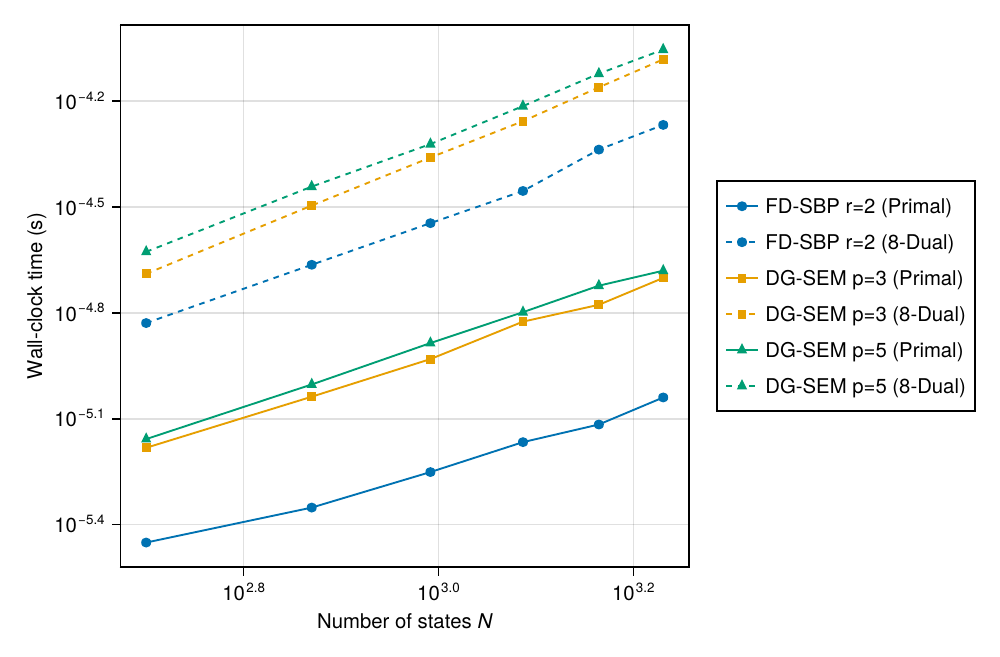}
\caption{Wall-clock time comparison for right-hand-side evaluation using primal and 8-partial dual state vectors.}
\label{fig:rhs_evaluation_time}
\end{figure}

Although the programmed inlet profile $c_{\mathrm{prog},i}(t)$ is discontinuous, it is not imposed directly on the column. First, the mixer in Eq.~\eqref{eq:cstr} smooths $c_{\mathrm{prog},i}(t)$ such that the concentration entering the column, $c_{\mathrm{in},i}(t)$, is continuous. Second, the Danckwerts inlet condition in Eq.~\eqref{eq:danckwerts_in} relates $c_{\mathrm{in},i}(t)$ to the interior solution through a convective--dispersive flux balance rather than prescribing it as a Dirichlet value, so finite axial dispersion smooths the inlet profile further. Consequently, the boundary-node value $c_i(t,0)$ is continuously differentiable in time, and some of the high-order convergence observed for the sufficiently smooth manufactured problems may be retained.

The simulated chromatogram is shown in Figure~\ref{fig:chromatogram}, where the OD of each component is plotted against time. The small initial salt bump observed is due to proteins exchanging salt ions during loading. 

\subsubsection{Per-operation benchmarks of the ODE solver} \label{sec:per_op}

In this section, the per-operation components of ODE integration are benchmarked, including right-hand-side evaluation, sparse state Jacobian evaluation and linear solver performance for both primal and dual-valued runs. The dual-valued solve uses double precision with eight partial derivatives. These operations constitute the main contributors to the total simulation cost when integrating ODEs with \texttt{DifferentialEquations.jl} and are generic to any method-of-lines chromatography model, whereas their relative weights depend on the operating regime. For a specific chromatography problem, its operating regime will determine how these individual costs combine to form the overall computational cost of simulation.

The operators are benchmarked for a range of grid sizes $N_i \in \{ 24, 36, 48, 60, 72, 84 \}$. DG-SEM fixes the polynomial degree $p \in \{3, 5\}$ and varies the number of elements, while FD-SBP fixes the boundary order $r=2$ and varies the number of nodes. With $N_c=6$ proteins, the column contributes $3(N_c+1)N_i - N_i$ states, since the salt loading $q_s$ is not carried as a state, and the mixer contributes $N_c+1$ states. Together with $2(N_c+1)$ additional states reserved for mass collection, the total number of states is $N = 20N_i + 21$, i.e.\ $N \in \{501, 741, 981, 1221, 1461, 1701\}$.

\begin{figure}[t!]
\centering
\includegraphics[width=0.90\textwidth]{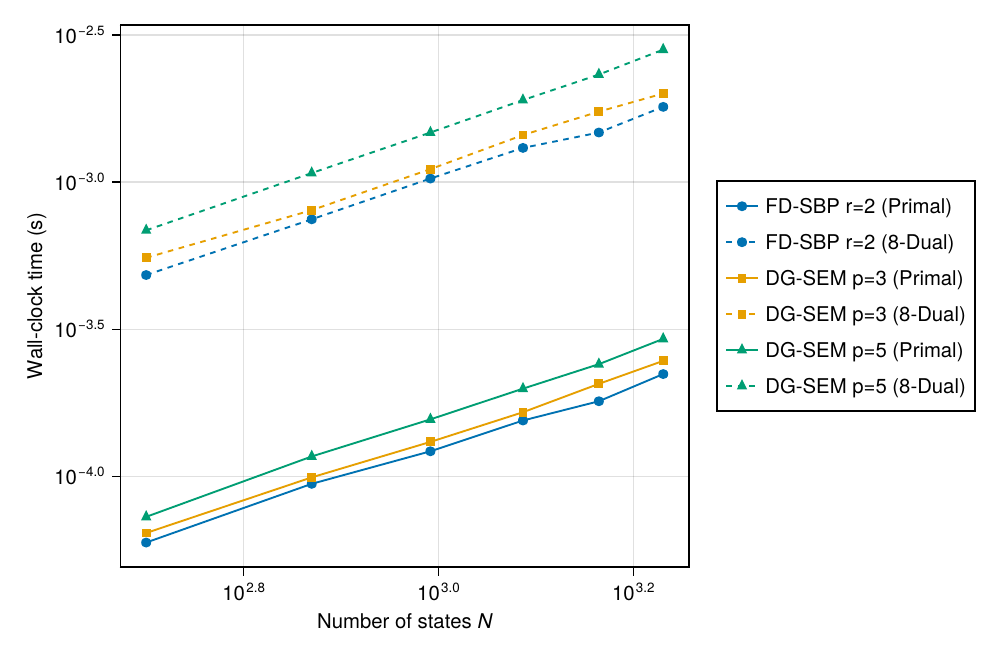}
\caption{Wall-clock time comparison for sparse Jacobian evaluation using primal and 8-partial dual state vectors.}
\label{fig:jacobian_evaluation_time}
\end{figure}

\paragraph{Semi-discrete right-hand-side evaluation:}

The semi-discrete right-hand-side $\dot{\mathbf{y}} = \mathbf{f}(\mathbf{y}, t)$ is implemented without memory allocations for performance. Since it is evaluated repeatedly during time integration, it is a natural first component to isolate.

Figure~\ref{fig:rhs_evaluation_time} reports the wall-clock time for one primal and one 8-partial dual-valued right-hand-side $\mathbf{f}(\mathbf{y}, t)$ evaluation. FD-SBP is faster than DG-SEM at every $N$ in both the primal and dual-valued evaluations. The cost of evaluation is roughly similar for DG-SEM for $p=3$ and $p=5$. The dual/primal overhead rises from $4.0$ at $N=501$ to $5.5$ at $N=1701$ for FD-SBP but stays in a tighter $3.2$--$4.0$ band for DG-SEM at the same $N$ and $p=3$. Both overheads are well below the naive factor of $9$ expected from carrying eight partials alongside the primal value, indicating that SIMD vectorization is effective for both operators. DG-SEM's tighter overhead is consistent with an element-local working set that does not grow with $N$, whereas FD-SBP's working set grows with $N$.

\paragraph{Sparse Jacobian evaluation:}

The sparse Jacobian is assembled with coloured forward-mode automatic differentiation using a greedy graph-colouring algorithm~\citep{Montoison2025SparseMatrixColorings}. Due to element local operations, the DG-SEM structural pattern for the mobile-phase self-interaction block $c_i \to c_i$ is block-diagonal with dense blocks corresponding to each element. In contrast, the mobile-phase self-interaction block $c_i \to c_i$ for the FD-SBP structural pattern is a narrow band due to the bandwidth optimized FD-SBP operators used. The remaining structural pattern entries are similar for the two methods. 

Figure~\ref{fig:jacobian_evaluation_time} reports the wall-clock time of a single sparse Jacobian $J$ evaluation, for both primal and 8-partial dual state vectors. FD-SBP is slightly faster than DG-SEM using $p=3$, but the difference is modest. For DG-SEM using $p=5$, the evaluation cost increases, due to a denser block-diagonal structure in the mobile-phase self-interaction block. Hence, for this operator to remain efficient, the problem must be sufficiently smooth to justify the higher polynomial degree.

\begin{figure}[t!]
\centering
\includegraphics[width=0.90\textwidth]{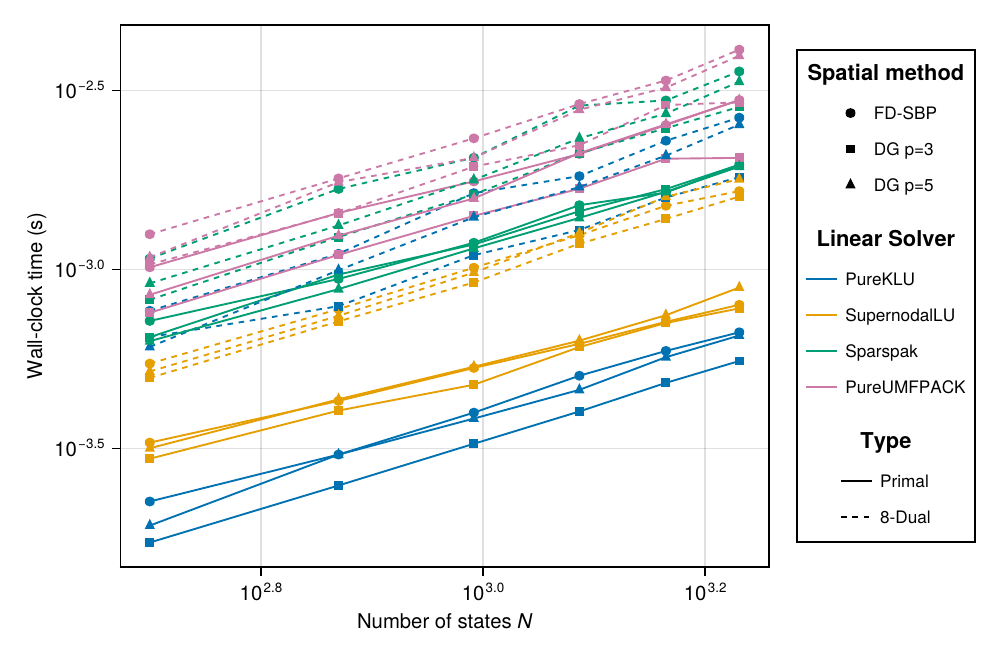}
\caption{Wall-clock time comparison for numerical refactorization using primal and 8-partial dual state vectors.}
\label{fig:linsolve_numerical_refactorization}
\end{figure}

\paragraph{Linear solver benchmark:}

The sparse linear solve benchmark is divided into its two operations, numerical refactorization and triangular solve. The total linear algebra cost is a weighted sum of these two components, depending on how often the chosen ODE integrator needs to refactorize versus performing triangular solves.

Four sparse direct solvers are benchmarked through the \texttt{LinearSolve.jl} interface, all implemented in Julia: \texttt{PureKLU.jl}, a bit-for-bit Julia port of KLU \citep{Davis2010KLU}; \texttt{PureUMFPACK.jl}, a Julia translation of UMFPACK~\citep{Davis2004UMFPACK}; \texttt{Sparspak.jl}, a sparse LU with minimum degree ordering \citep{GeorgeLiu1981Sparspak}; and \texttt{SupernodalLU}, a left-looking supernodal LU \citep{SCHENK2004475, Schenk2006}, vendored inside \texttt{LinearSolve.jl}. BLAS/LAPACK \citep{Dongarra1990Level3BLAS} calls use the Intel MKL backend via \texttt{MKL.jl}.

\begin{figure}[t!]
\centering
\includegraphics[width=0.90\textwidth]{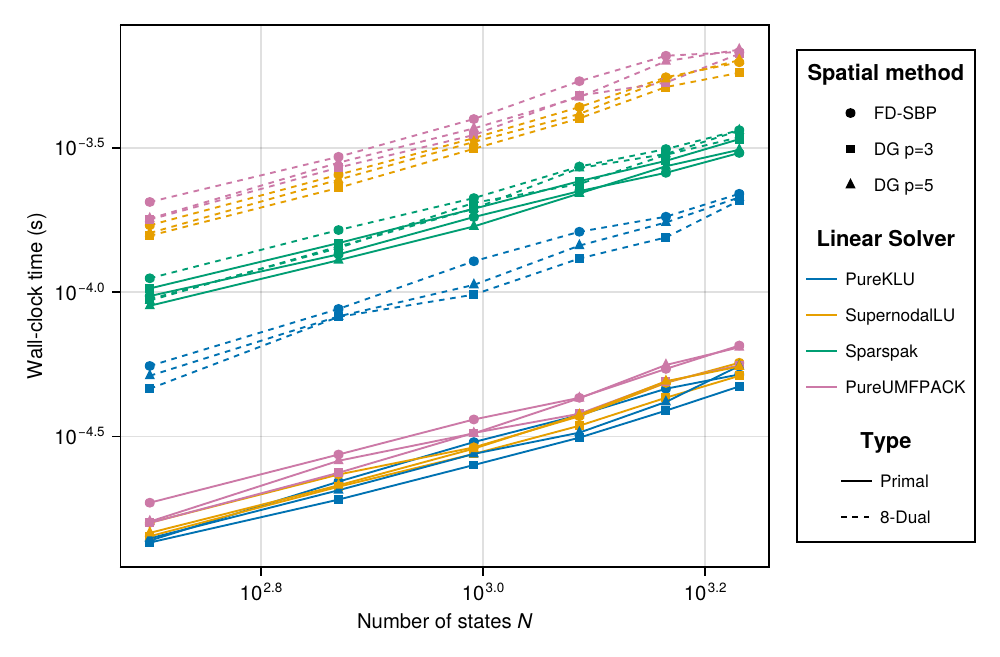}
\caption{Wall-clock time comparison for triangular solve using primal and 8-partial dual state vectors.}
\label{fig:linsolve_triangular_solve}
\end{figure}

Figure~\ref{fig:linsolve_numerical_refactorization} and Figure~\ref{fig:linsolve_triangular_solve} shows the cost of primal and dual-valued refactorization and triangular solve, respectively, for the different sparse direct solvers. For numerical refactorization, \texttt{PureKLU.jl} performs best with primal values, while \texttt{SupernodalLU} performs best with dual values. For triangular solve, \texttt{PureKLU.jl}, \texttt{SupernodalLU} and \texttt{PureUMFPACK.jl} are closely matched with primal values, while \texttt{PureKLU.jl} performs best with dual values. Hence, for primal simulations, \texttt{PureKLU.jl} is clearly the most efficient choice. For dual-valued simulations, \texttt{SupernodalLU} is slightly faster in the numerical refactorization, but \texttt{PureKLU.jl} is faster in the triangular solve. Since the numerical refactorization is typically performed less frequently than the triangular solve, \texttt{PureKLU.jl} is also the preferred choice for dual-valued simulations. In the remainder of the paper, \texttt{PureKLU.jl} will be used as the default sparse direct solver for both primal and dual-valued simulations.

In both refactorization and triangular solve, DG-SEM with $p=3$ has the lowest cost followed by DG-SEM with $p=5$ and finally FD-SBP is slightly more expensive with $r=2$. The main difference arises from the sparsity pattern of the Jacobian, where DG-SEM benefits from a block-diagonal structure with dense blocks, while FD-SBP has a narrow band structure that slightly increases the computational cost. As the polynomial degree is increased in DG-SEM, the blocks become more dense, leading to increased computational cost. As for Jacobian evaluation, the problem must be sufficiently smooth to justify the higher polynomial degree.

\subsubsection{Primal solve benchmark} \label{sec:primal_solve_benchmark}

\begin{table}[t!]
\centering
\caption{Cost breakdown of the primal QNDF work-precision sweep at matched $N$ (DG-SEM vs FD-SBP). Components: right-hand-side (RHS), Jacobian (Jac), refactorization (Refac), triangular solve (Tri), and Other. The FD-SBP interior order is $l=2r$. Figure~\ref{fig:wp_dg_vs_fdsbp_time} shows the total wall-clock time.}
\label{tab:wp_cost_decomp}
\begin{tabular}{lrrrrrrrrr}
\toprule
Method & $p$ & $r$ & $N$ & RMS-error & RHS  & Jac  & Refac  & Tri  & Other  \\
       &     &     &     & (\unit{\AU\per\centi\meter})   & (\%) & (\%) & (\%)   & (\%) & (\%)   \\
\midrule
DG-SEM & 5 & --- &  501 & 1.614\text{e--}2 & 22.6 & 4.0 & 15.6 & 39.9 & 17.9 \\
DG-SEM & 5 & --- &  741 & 7.649\text{e--}3 & 20.2 & 4.2 & 15.1 & 39.9 & 20.7 \\
DG-SEM & 5 & --- &  981 & 3.801\text{e--}3 & 21.8 & 3.8 & 14.2 & 43.6 & 16.6 \\
DG-SEM & 5 & --- & 1221 & 1.738\text{e--}3 & 20.8 & 4.0 & 14.5 & 41.5 & 19.1 \\
\addlinespace
FD-SBP & --- & 2 &  501 & 8.509\text{e--}3 & 12.0 & 3.0 & 19.1 & 44.0 & 21.9 \\
FD-SBP & --- & 2 &  741 & 2.233\text{e--}3 & 11.7 & 3.3 & 18.3 & 50.0 & 16.6 \\
FD-SBP & --- & 2 &  981 & 7.422\text{e--}4 & 10.4 & 3.5 & 17.4 & 46.6 & 22.1 \\
FD-SBP & --- & 2 & 1221 & 3.186\text{e--}4 &  9.8 & 3.6 & 19.0 & 54.1 & 13.5 \\
\bottomrule
\end{tabular}
\end{table}

At matched $N_i$, the three primal per-operation benchmarks show that FD-SBP evaluates the right-hand-side and Jacobian faster than DG-SEM, whereas DG-SEM is faster for the sparse linear solve. The two spatial discretizations share key mathematical properties -- both are energy-stable and high-order accurate -- but differ in implementation. They can therefore be viewed as complementary approaches with distinct computational trade-offs across these operations. Their combined effect on full-solve performance is assessed next. 

Since no analytical solution is available, the spatial discretization
error is quantified against a highly resolved reference solution $c_{\mathrm{ref},i}(t,L)$
obtained with DG-SEM on 400 elements of degree $p=4$ with \texttt{QNDF} at
$\mathrm{abstol} = \mathrm{reltol} = 10^{-12}$. Accuracy is measured by the
root-mean-square (RMS) error of the OD outlet signal
over $[0,t_f]$,
\begin{equation}
\mathrm{RMS\text{-}error} = w\sqrt{\frac{1}{t_f}\int_0^{t_f}
\sum_{i=1}^{N_c} \big(c_i(t,L)-c_{\mathrm{ref},i}(t,L)\big)^2\,dt}.
\label{eq:rms_metric}
\end{equation}
The time integral is evaluated by sampling the outlet signals at $1000$ time points uniformly spaced over $[0,t_f]$, fitting a shape-preserving PCHIP interpolant to the pointwise squared error with \texttt{DataInterpolations.jl}, and integrating it analytically.

\begin{figure}[t!]
\centering
\includegraphics[width=0.90\textwidth]{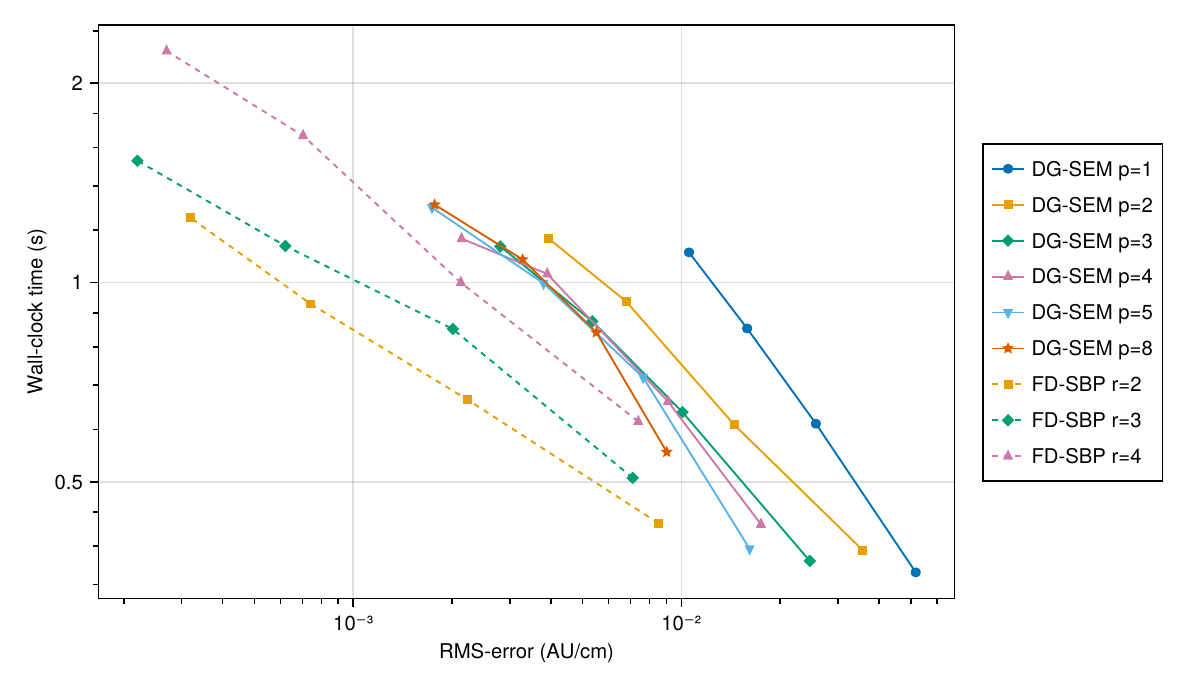}
\caption{Primal work-precision comparison between DG-SEM and FD-SBP implementations. See details in Table~\ref{tab:wp_cost_decomp}.}
\label{fig:wp_dg_vs_fdsbp_time}
\end{figure}

Three checks verify that the reference solution is sufficiently resolved. First, recomputing the solution on a coarser, 200-element mesh changes the RMS error by only $7.6\times10^{-8}$~\unit{\AU\per\centi\meter} relative to the 400-element mesh, indicating that the spatial discretization is effectively converged; the finer 400-element solution is therefore retained as the reference. Second, recomputing the solution on the 200-element grid with a tenfold tighter time tolerance changes the RMS error by only $4.2\times10^{-8}$~\unit{\AU\per\centi\meter}, indicating that temporal error is negligible. Finally, a highly resolved FD-SBP solution ($r=2$, 600 nodes) agrees with the 400-element DG-SEM reference to within $7.4\times10^{-7}$~\unit{\AU\per\centi\meter}, confirming method neutrality.

The work-precision sweep is performed by varying the number of grid nodes $N_i \in \{24, 36, 48, 60 \}$ for both DG-SEM and FD-SBP. The DG-SEM sweep uses polynomial degrees $p \in \{ 1, 2, 3, 4, 5, 8 \}$, while the FD-SBP sweep uses boundary orders $r \in \{2, 3, 4 \}$. The degrees $p=4$ and $p=8$ do not divide these target sizes evenly, so the DG-SEM sweep instead uses the nearest admissible resolutions $N_i \in \{25, 35, 50, 60\}$ and $N_i \in \{27, 36, 45, 54\}$, respectively. The QNDF integrator is used with $\mathrm{reltol}=\mathrm{abstol}=10^{-10}$, sufficient to make temporal errors negligible.

The resulting RMS-errors and wall-clock times are reported in Figure~\ref{fig:wp_dg_vs_fdsbp_time}, showing that FD-SBP consistently outperforms DG-SEM on this benchmark. FD-SBP is most efficient at boundary order $r=2$, with diminishing returns from further increases in order. DG-SEM shows a similar trend as the polynomial degree $p$ increases, with only marginal gains beyond $p=3$; $p=5$ is the most accurate at a given $N_i$ and is therefore used in the following. These results suggest that high-order convergence is retained, but that the limited smoothness of the solution restricts the observed convergence rate. The pairwise EOCs computed from Table~\ref{tab:wp_cost_decomp} are $3.3$--$3.8$ for FD-SBP ($r=2$) and increase from $1.8$ to $3.5$ for DG-SEM ($p=5$). This is well below the design order of $p+1$ for DG-SEM. Increasing the order further therefore offers little advantage on this benchmark.

For a similar grid $N_i$, the FD-SBP operator reaches a lower RMS-error than DG-SEM. This advantage may be attributed in part to its continuous stencil, which avoids the duplication of DOF at element interfaces as in DG-SEM. Unless the DG-SEM formulation can exploit the discontinuities between elements, these additional DOF primarily increase the computational cost without a corresponding gain in accuracy. 

Table~\ref{tab:wp_cost_decomp} decomposes each swept case into contributions from the right-hand side, sparse Jacobian evaluation, numerical refactorization and triangular solve. Each unit cost is estimated by the per-operation time multiplied by the number of times the operation is invoked during the QNDF solve. The remainder of the total wall-clock time collects interpolation, output storage, etc. 

Table~\ref{tab:wp_cost_decomp} shows that right-hand-side (RHS) evaluation accounts for $10$--$12$\% of the total cost for FD-SBP and $20$--$23$\% for DG-SEM. Jacobian evaluation contributes only $3$--$4$\% for both methods. QNDF reuses matrix factorizations, limiting numerical refactorization to $17$--$19$\% for FD-SBP and $14$--$16$\% for DG-SEM. Triangular solves dominate the linear algebra cost, accounting for $44$--$54$\% for FD-SBP and $40$--$44$\% for DG-SEM. Thus, RHS evaluations account for a smaller share of runtime in FD-SBP, while linear algebra accounts for a smaller share in DG-SEM, as expected from the per-operation benchmarks (section~\ref{sec:per_op}).

Jacobian evaluation could be accelerated by precomputing the constant convection and dispersion contributions and updating only the entries associated with the nonlinear SMA isotherm \citep{Meyer2020ChromaTech}. However, since the cost of evaluating the Jacobian is almost negligible, the potential speed-up from this optimization is limited. For this reason, automatic differentiation is retained for simplicity.

The results above are obtained for a single, representative operating point. Some of the observed differences are structural and are expected to carry over to other parameter values: the lower number of DOF of FD-SBP for a given resolution and its cheaper right-hand-side evaluation do not depend on the model parameters. In contrast, the relative weights of the individual costs depend on the operating regime. When linear algebra dominates, e.g.\ for stiffer adsorption kinetics or a larger number of components, the linear-algebra advantage of DG-SEM becomes more important. Moreover, for strongly convection-dominated problems (large $\mathrm{Pe}$) on under-resolved grids, the dissipation introduced by the upwind numerical flux of DG-SEM may be advantageous.

\begin{figure}[t!]
\centering
\includegraphics[width=0.90\textwidth]{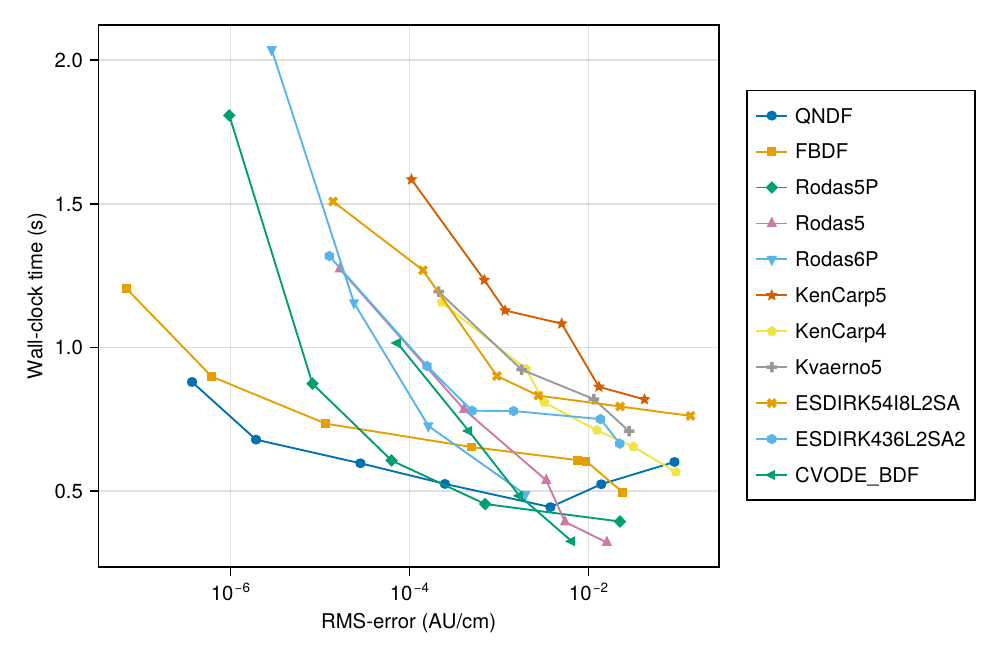}
\caption{Work-precision benchmark for ODE integrators using a FD-SBP semi-discrete scheme with boundary order $r=2$ and $40$ nodes.}
\label{fig:wp_high_order_l2}
\end{figure}

\subsubsection{ODE time integration benchmark}

Having considered the spatial discretization, the remaining algorithmic choice is the ODE integrator. This benchmark uses a FD-SBP semi-discrete scheme with boundary order $r=2$ and $40$ nodes and compares how different time integrators perform. Because both the test and reference solutions are computed on the same spatial grid, the measured error reflects time-integration error alone. In the multi-step BDF family, the Julia solvers QNDF and FBDF are tested and the CVODE\_BDF method from Sundials \citep{Hindmarsh2005SUNDIALS}. In the DIRK family, Kvaerno5, KenCarp4/5, ESDIRK54I8L2SA and ESDIRK436L2SA2 are tested. In the Rosenbrock--Wanner family, Rodas5, Rodas5P and Rodas6P are tested.

The reference solution is obtained by integrating the semi-discrete system with \texttt{Rodas6P} at $\mathrm{abstol}=\mathrm{reltol}=10^{-12}$. Tightening the tolerances to $\mathrm{abstol}=\mathrm{reltol}=10^{-13}$ changes the RMS-error by only $1.5\times10^{-10}$ \unit{\AU\per\centi\meter}, confirming that the reference is temporally resolved.

\begin{figure}[t!]
\centering
\includegraphics[width=0.90\textwidth]{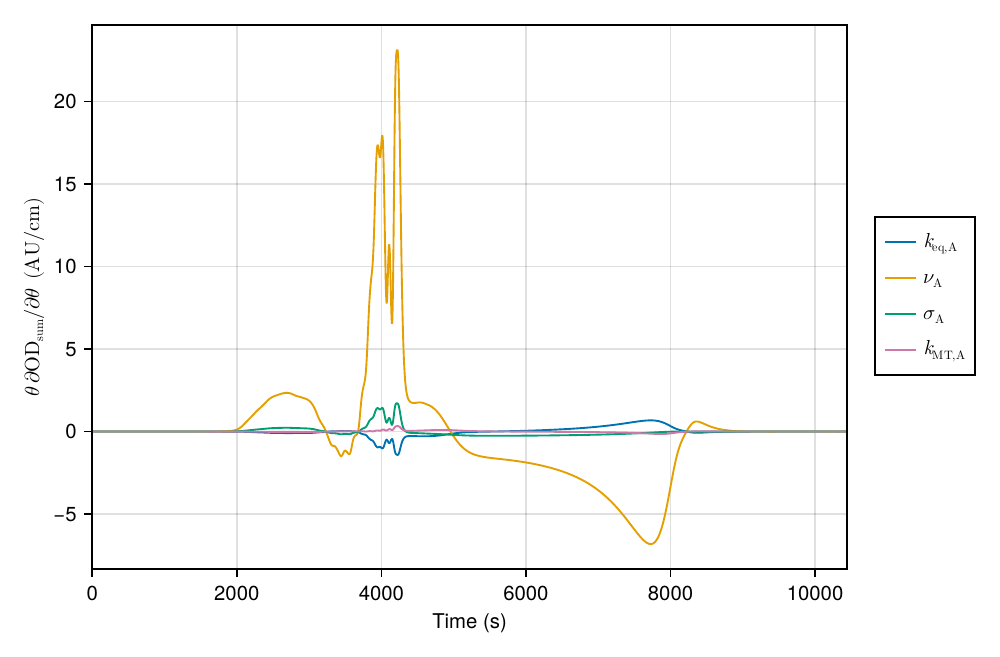}
\caption{Scaled sensitivities of the summed OD signal to the four parameters of component A for the 6-component ion-exchange problem. The corresponding chromatogram is shown in Figure~\ref{fig:chromatogram}.}
\label{fig:sensitivities}
\end{figure}

Figure~\ref{fig:wp_high_order_l2} shows work-precision performance for the three integrator families. DIRK methods are not competitive over the tested tolerance range, whereas Rodas5P remains competitive at intermediate tolerances. Tighter tolerances favor the BDF methods, particularly QNDF. This advantage may reflect more effective reuse of matrix factorizations across time steps, reducing the computational cost relative to Rosenbrock--Wanner methods. This balance may shift toward Rosenbrock--Wanner methods when QNDF requires frequent refactorizations or costly Newton iterations, reducing the benefits of factorization reuse across time steps. In the remainder of the paper, QNDF is used as the time integrator.

\subsubsection{Discrete forward sensitivity benchmark}

\begin{figure}[t!]
\centering
\includegraphics[width=0.90\textwidth]{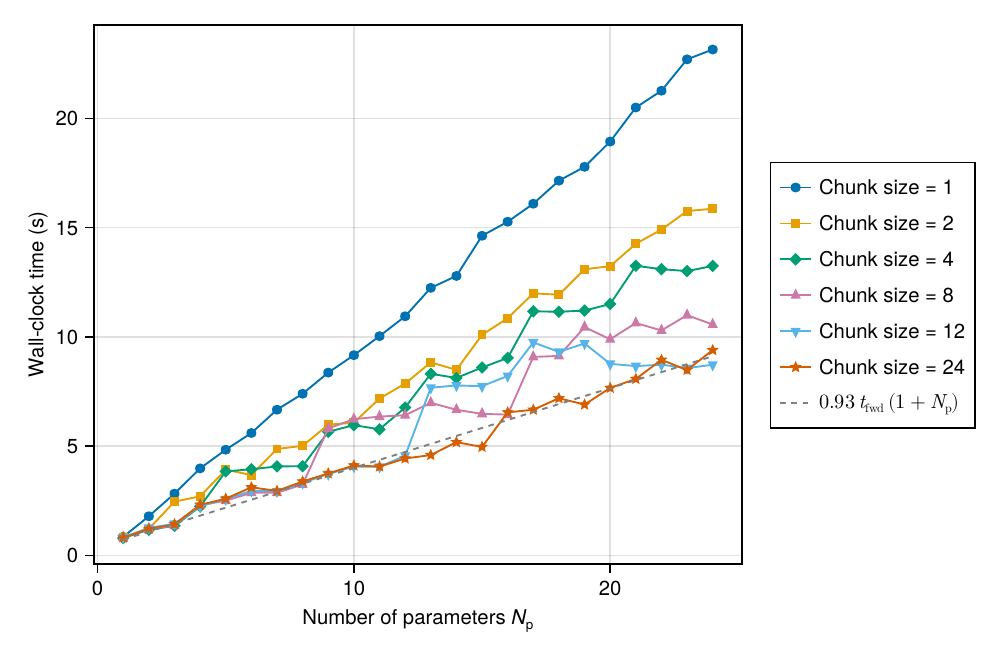}
\caption{DFSA cost vs propagated parameters $N_p$ for fixed \texttt{ForwardDiff.jl} chunk widths $N_\mathrm{chunk}\in\{1,2,4,8,12,24\}$. Dashed line: single chunk estimator $k\,t_{\mathrm{fwd}}(1+N_p)$, with $t_{\mathrm{fwd}}$ the forward simulation time.}
\label{fig:dsa_scaling_chunks}
\end{figure}

The dual per-operation benchmarks mirror the primal case: FD-SBP is faster on RHS and Jacobian evaluations whereas DG-SEM is faster on linear algebra. Whether these effects combine to alter the ranking in a full dual solve is assessed next. In this section, the ODE integrator and linear solver are fixed to \texttt{QNDF} and \texttt{PureKLU.jl}, respectively

The parameters are the SMA isotherm coefficients ($k_{\mathrm{eq}, i}$, $\nu_i$, $\sigma_i$) and the mass-transfer coefficient ($k_{\mathrm{MT}, i}$). They are ordered componentwise as
\[
\big(k_{\mathrm{eq},\mathrm{A}},\, \nu_{\mathrm{A}},\, \sigma_{\mathrm{A}},\, k_{\mathrm{MT},\mathrm{A}},\, \ldots,\,
k_{\mathrm{eq},\mathrm{F}},\, \nu_{\mathrm{F}},\, \sigma_{\mathrm{F}},\, k_{\mathrm{MT},\mathrm{F}}\big),
\]
so that a benchmark with $N_p$ parameters uses the first $N_p$ entries of this list, up to a maximum of 24. Here, $k_{\mathrm{MT},i}$ is treated as a constant. If it instead depends on the flow rate through a correlation $k_{\mathrm{MT},i}(v_{\mathrm{int}};\boldsymbol{\theta})$, the correlation simply becomes part of the differentiated program: the sensitivity with respect to $v_{\mathrm{int}}$ then automatically includes the indirect contribution through $k_{\mathrm{MT},i}$, and the correlation coefficients $\boldsymbol{\theta}$ replace $k_{\mathrm{MT},i}$ as parameters.

The output is the summed OD signal $\text{OD}_\text{sum}(t; \boldsymbol{\theta}) = \sum_{i=1}^{N_c} w\,c_i(t,L; \boldsymbol{\theta})$, and the benchmark measures the wall-clock time required to compute the associated $N_p$ parameter sensitivities $\partial \text{OD}_\text{sum}(t; \boldsymbol{\theta}) / \partial \theta_j$ for $j=1,\ldots,N_p$. This illustrates the ease of applying DFSA: the user simply defines the desired output as a function of the state vector, and \texttt{ForwardDiff.jl} automatically propagates the sensitivities through the same code path. Changing the output, the parameters, or the model equations themselves therefore requires no hand-derived sensitivity code, only a re-run of the differentiated program. 

Figure~\ref{fig:sensitivities} shows the scaled sensitivities $S_k(t) = \partial \text{OD}_\text{sum}(t)/\partial\theta_k\,\theta_k$ for the four parameters of component A ($k_{\mathrm{eq},\mathrm{A}}$, $\nu_\mathrm{A}$, $\sigma_\mathrm{A}$, $k_{\mathrm{MT},\mathrm{A}}$). This scaling allows the effects of equal percentage changes in different parameters to be compared directly. For example, increasing one parameter by $1$\%, while holding the others fixed, changes the outlet signal by approximately $0.01S_k(t)$ for sufficiently small parameter changes. The magnitude $|S_k(t)|$ shows how strongly the signal responds to a small relative change in $\theta_k$, while its sign shows whether the signal increases or decreases. Sensitivities are largest near steep elution fronts, where a small shift in elution time causes a large change in $\text{OD}_\text{sum}$. Since the injected proteins are completely eluted within the simulated time, the total peak area is conserved and each sensitivity integrates to zero. A positive region followed by a negative region is therefore consistent with earlier elution since the signal increases on the leading front and decreases on the tail. For instance, increasing $\nu_\mathrm{A}$ moves the elution of A earlier, raising the signal on its leading front and lowering it on its tail. Since the proteins compete for binding sites, the parameters of A also shift the other components, explaining the nonzero sensitivities outside A's elution window.

\paragraph{Chunk-width scaling:} The first benchmark examines how the \texttt{ForwardDiff.jl} chunk width affects the cost of computing the $N_p$ parameter sensitivities. To isolate this effect, the spatial discretization is fixed to FD-SBP with boundary order $r=2$ and $40$ nodes. The time tolerances are set to $\mathrm{abstol}=\mathrm{reltol}=10^{-7}$.

The parameter sensitivities are computed as the parameter count $N_p$ is swept from $1$ to $24$, with the chunk width pinned to each of $N_\mathrm{chunk}\in\{1,2,4,8,12,24\}$ (clamped to $N_p$ for narrow problems). These widths are chosen so that they divide $24$ evenly. For 256-bit SIMD with \texttt{Float64} values (4 doubles per vector), $N_\mathrm{chunk}\in\{4,8,12,24\}$ maps to an integer number of full vectors and avoids padding, whereas $N_\mathrm{chunk}\in\{1,2\}$ does not.

Figure~\ref{fig:dsa_scaling_chunks} reports the resulting wall-clock times. With a fixed chunk width, the gradient is assembled over $\lceil N_p/N_\mathrm{chunk}\rceil$ dual passes, each re-running the primal solve while propagating up to $N_\mathrm{chunk}$ partials. The number of passes, not $N_p$ itself, is therefore the main cost driver: adding a parameter that still fits within the current passes costs only one extra partial, whereas adding one that starts a new pass costs a full extra chunk including unused padding. This explains the staircase pattern in the computational cost curves. 

The single-chunk configuration ($N_\mathrm{chunk}=N_p$, one pass) is generally cheapest, with two exceptions where a narrower chunk is competitive: at $N_p=16$, two passes with $N_\mathrm{chunk}=8$ is competitive, and at $N_p=24$, two passes with $N_\mathrm{chunk}=12$ is competitive. In these regimes the wide single chunk incurs enough per-partial overhead (e.g., larger dual objects and less favourable cache behaviour) that multiple passes with zero padding can be competitive, despite having to compute the primal solve twice. The dashed line in Figure~\ref{fig:dsa_scaling_chunks} shows the single-chunk estimator $0.93\,t_{\mathrm{fwd}}\,(1+N_p)$, i.e.\ about $0.93$ forward solves per parameter. This suggests that the chunked approach is effective since it can enable SIMD vectorization.

Across the entire sweep, any chunk width larger than $1$ outperforms $N_\mathrm{chunk}=1$ by a wide margin, so using $N_\mathrm{chunk}>1$ is always preferable. The width $N_\mathrm{chunk}=2$ also rarely pays off. It is only competitive with the $N_\mathrm{chunk}=4$ and $N_\mathrm{chunk}=8$ configurations where they carry the most unnecessary padding. 

\paragraph{Work-precision benchmark:}

\begin{figure}[t!]
\centering
\includegraphics[width=0.90\textwidth]{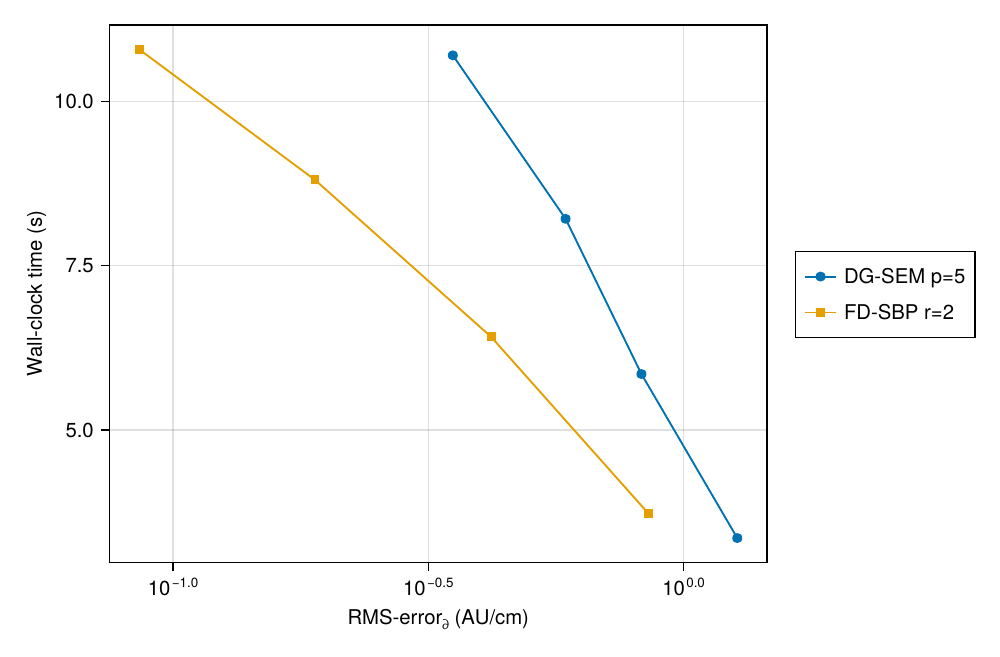}
\caption{Dual-valued work-precision comparison between DG-SEM and FD-SBP implementations. $N_p=N_\mathrm{chunk}=8$. See details in Table~\ref{tab:wp_dsa_cost_decomp}. See the primal counterpart in Fig.~\ref{fig:wp_dg_vs_fdsbp_time}.}
\label{fig:wp_dsa_l2}
\end{figure}

\begin{table}[t!]
\centering
\caption{Cost breakdown of the 8-partial dual-valued QNDF work-precision sweep at matched $N$ (DG-SEM vs FD-SBP) with $\mathrm{abstol}=\mathrm{reltol}=10^{-10}$. Components: right-hand-side (RHS), Jacobian (Jac), refactorization (Refac), triangular solve (Tri), solver setup (Setup), and Other. The FD-SBP interior order is $l=2r$. Figure~\ref{fig:wp_dsa_l2} shows the total wall-clock time. Table~\ref{tab:wp_cost_decomp} shows the corresponding breakdown for the primal work-precision sweep.}
\label{tab:wp_dsa_cost_decomp}
\begin{tabular}{lrrrrrrrrrr}
\toprule
Method & $p$ & $r$ & $N$ & RMS-error  & RHS  & Jac  & Refac  & Tri  & Setup & Other  \\
       &     &     &     & (\unit{\AU\per\centi\meter})  & (\%) & (\%) & (\%)   & (\%) & (\%)  & (\%)   \\
\midrule
DG-SEM & 5 & --- &  501 & 1.27\text{e+}0  & 11.0 & 4.5 & 8.4 & 18.4 & 0.3 & 57.5 \\
DG-SEM & 5 & --- &  741 & 8.27\text{e--}1 & 10.6 & 4.3 & 7.1 & 21.1 & 0.2 & 56.7 \\
DG-SEM & 5 & --- &  981 & 5.87\text{e--}1 & 10.5 & 4.6 & 7.5 & 18.9 & 0.1 & 58.3 \\
DG-SEM & 5 & --- & 1221 & 3.53\text{e--}1 & 11.2 & 4.5 & 7.4 & 20.5 & 0.1 & 56.3 \\
\addlinespace
FD-SBP & --- & 2 &  501 & 8.52\text{e--}1 & 10.0 & 3.2 & 9.3 & 19.9 & 0.2 & 57.5 \\
FD-SBP & --- & 2 &  741 & 4.20\text{e--}1 &  8.3 & 3.3 & 8.1 & 22.1 & 0.1 & 58.2 \\
FD-SBP & --- & 2 &  981 & 1.90\text{e--}1 &  9.3 & 3.1 & 8.3 & 21.5 & 0.1 & 57.9 \\
FD-SBP & --- & 2 & 1221 & 8.59\text{e--}2 &  8.3 & 3.4 & 9.2 & 24.6 & 0.1 & 54.5 \\
\bottomrule
\end{tabular}
\end{table}

The chunk-width study quantifies the cost scaling of dual propagation for a fixed solver configuration, but it does not address the accuracy-to-cost trade-off. To evaluate this trade-off, FD-SBP and DG-SEM are compared on an $N_p=8$ parameter problem.

To quantify accuracy, the OD-weighted RMS metric of the sensitivities is defined analogously to Eq.~\eqref{eq:rms_metric} as
\begin{equation}
\mathrm{RMS\text{-}error}_{\partial} = \sqrt{\frac{1}{t_f}\int_0^{t_f}
\sum_{k=1}^{N_p} \left(\theta_k\left(
\frac{\partial \textrm{OD}_\textrm{sum}(t; \boldsymbol{\theta})}{\partial\theta_k}
- \frac{\partial \textrm{OD}_\textrm{sum}^{\textrm{ref}}(t; \boldsymbol{\theta})}{\partial\theta_k}\right)\right)^2 dt}.
\end{equation}

The work-precision sweep parameters are chosen to be equal to those of the primal sweep described in Section~\ref{sec:primal_solve_benchmark}, using the same grid node configurations and time tolerances. Based on the primal solve, only the best performing boundary order $r=2$ and polynomial degree $p=5$ are investigated. The reference solution is confirmed adequate using the same procedure as for the primal solve.

Figure~\ref{fig:wp_dsa_l2} shows the resulting work-precision diagram. Consistent with the primal sweep, FD-SBP outperforms DG-SEM. Table~\ref{tab:wp_dsa_cost_decomp} shows that dual-valued RHS evaluation accounts for $8$--$10$\% of the total cost for FD-SBP and $10$--$11$\% for DG-SEM. Jacobian evaluation contributes only $3$\% for FD-SBP and $4$--$5$\% for DG-SEM. Numerical refactorization accounts for $8$--$9$\% for FD-SBP and $7$--$8$\% for DG-SEM. Again, triangular solves dominate the linear algebra cost, accounting for $20$--$25$\% for FD-SBP and $18$--$21$\% for DG-SEM. The Other component, which collects all operations not captured by the timed kernels, e.g.\ storing the dual-valued output, is significantly larger for the dual-valued solve than for the primal solve.

\section{Conclusion}\label{sec:conclusion}

This paper presented \texttt{ChromOps.jl}, a Julia-based Novo Nordisk A/S inner-source chromatography simulation framework that combines high-order spatial discretization
with DFSA. The framework was benchmarked on a preparative 6-component ion-exchange
problem using a lumped-rate model with SMA isotherm kinetics.

Comparing spatial discretizations, FD-SBP outperforms DG-SEM both in primal and 8-partial dual-valued solves on the considered benchmark. Both attain their theoretical high-order convergence rates on manufactured solutions with sufficient smoothness. The per-operation computational costs, including right-hand-side and sparse Jacobian evaluation, refactorization and triangular solve, show that FD-SBP is faster on the right-hand-side and Jacobian evaluation, whereas DG-SEM is faster on linear algebra for similar grid sizes $N_i$. The FD-SBP method achieves higher accuracy per DOF than DG-SEM, which can be attributed in part to its continuous stencil that avoids duplicated interface nodes. Unless the inter-element jumps of DG-SEM are exploited, e.g.\ through flux dissipation in under-resolved simulations, the duplicated interface DOF are a disadvantage in terms of work-precision performance. Moreover, the derivation and implementation of FD-SBP require substantially less effort than DG-SEM. For these reasons, it offers a competitive alternative to DG-SEM for chromatography simulation.

The DFSA was benchmarked for up to $N_p=24$ parameters, with the cost of propagating $N_p$ parameters scaling roughly linearly in $N_p$. The chunk width of the \texttt{ForwardDiff.jl} propagation can be tuned to balance memory use, with a single wide chunk optimal in most of the sweep. The exceptions where multiple passes become competitive are for $N_p=16$ and $N_p=24$, where $N_\mathrm{chunk}=8$ and $N_\mathrm{chunk}=12$ splits the work into two full passes without padding.

Based on the benchmarks, the following settings are recommended as a starting point for similar chromatography problems: FD-SBP with boundary order $r=2$ and a sufficient number of nodes, e.g.\ $N_i\approx40$--$60$; the \texttt{QNDF} integrator with \texttt{PureKLU.jl} as sparse linear solver; and a single chunk $N_\mathrm{chunk}=N_p$, or an even split into two chunks when this avoids padding, e.g.\ $N_\mathrm{chunk}=8$ for $N_p=16$ and $N_\mathrm{chunk}=12$ for $N_p=24$. Since the sensitivities can converge more slowly than the primal solution under grid refinement (cf.\ Tables~\ref{tab:wp_cost_decomp} and~\ref{tab:wp_dsa_cost_decomp}), the spatial resolution should be verified for the sensitivities as well.

Because the primal and dual-valued solves share the same code path, DFSA yields fast and accurate gradients of user-defined objective functions without requiring the user to specify chain-rule code. Gradient-based parameter estimation and optimization thus become accessible to chromatography practitioners without expertise in sensitivity methods, supporting the model-based process development workflows at the core of chromatography PSE.

\section*{Funding}

This research did not receive any specific grant from funding agencies in the public, commercial, or not-for-profit sectors.

\section*{CRediT authorship contribution statement}

\textbf{Kristian Meyer:} Conceptualization, Methodology, Software, Validation, Formal analysis, Investigation, Data curation, Visualization, Writing -- original draft, Writing -- review \& editing. \textbf{Maksym Ratajczyk:} Conceptualization, Methodology, Software, Validation, Writing -- review \& editing. \textbf{Christopher Rackauckas:} Conceptualization, Software, Writing -- review \& editing.

\bibliographystyle{elsarticle-harv}
\bibliography{references}

@article{Kumar2020MechanisticReview,
  author  = {Kumar, Vijesh and Lenhoff, Abraham M.},
  title   = {Mechanistic Modeling of Preparative Column Chromatography for Biotherapeutics},
  journal = {Annu. Rev. Chem. Biomol. Eng.},
  year    = {2020},
  volume  = {11},
  pages   = {235--255},
  url     = {https://doi.org/10.1146/annurev-chembioeng-102419-125430}
}

@article{Benner2019LabToManufacturing,
  author  = {Benner, Steven W. and Welsh, John P. and Rauscher, Michael A. and Pollard, Jennifer M.},
  title   = {Prediction of lab and manufacturing scale chromatography performance using mini-columns and mechanistic modeling},
  journal = {J. Chromatogr. A},
  year    = {2019},
  volume  = {1593},
  pages = {54--62},
  url     = {https://doi.org/10.1016/j.chroma.2019.01.063}
}

@article{Rischawy2019GoodModelingPractice,
  author  = {Rischawy, Federico and Saleh, David and Hahn, Tobias and Oelmeier, Stefan and Spitz, Julia and Kluters, Simon},
  title   = {Good modeling practice for industrial chromatography: Mechanistic modeling of ion exchange chromatography of a bispecific antibody},
  journal = {Comput. Chem. Eng.},
  year    = {2019},
  volume  = {130},
  pages   = {106532},
  url     = {https://doi.org/10.1016/j.compchemeng.2019.106532}
}

@article{ranocha2022adaptive,
  title={Adaptive numerical simulations with {T}rixi.jl:
         {A} case study of {J}ulia for scientific computing},
  author={Ranocha, Hendrik and Schlottke-Lakemper, Michael and Winters, Andrew R.
          and Faulhaber, Erik and Chan, Jesse and Gassner, Gregor J.},
  journal={Proc. JuliaCon Conf.},
  volume={1},
  number={1},
  pages={77},
  year={2022},
  url={https://doi.org/10.21105/jcon.00077}
}

@article{schlottkelakemper2021purely,
  title={A purely hyperbolic discontinuous {G}alerkin approach for
         self-gravitating gas dynamics},
  author={Schlottke-Lakemper, Michael and Winters, Andrew R. and
          Ranocha, Hendrik and Gassner, Gregor J.},
  journal={J. Comput. Phys.},
  pages={110467},
  year={2021},
  volume={442},
  url={https://doi.org/10.1016/j.jcp.2021.110467}
}

@article{Ranocha2021SBPOperators,
  author    = {Ranocha, Hendrik},
  title     = {{SummationByPartsOperators.jl}: {A} {J}ulia library of provably stable discretization techniques with mimetic properties},
  journal   = {J. Open Source Softw.},
  year      = {2021},
  volume    = {6},
  number    = {64},
  pages     = {3454},
  url       = {https://doi.org/10.21105/joss.03454}
}

@article{Frandsen2025CADETJulia,
  title   = {{CADET-Julia}: Efficient and versatile, open-source simulator for batch chromatography in Julia},
  author  = {Frandsen, Jesper and Breuer, Jan M. and Schm{\"o}lder, Johannes and Huusom, Jakob K. and Gernaey, Krist V. and Abildskov, Jens and von Lieres, Eric},
  journal = {Comput. Chem. Eng.},
  year    = {2025},
  volume  = {192},
  pages   = {108913},
  url     = {https://doi.org/10.1016/j.compchemeng.2024.108913}
}

@article{Morbidelli1982POR,
  author  = {Morbidelli, Massimo and Servida, Alberto and Storti, Giuseppe and Carra, Sergio},
  title   = {Simulation of multicomponent adsorption beds. Model analysis and numerical solution},
  journal = {Ind. Eng. Chem. Fundam.},
  year    = {1982},
  volume  = {21},
  number  = {2},
  pages   = {123--131},
  url     = {https://doi.org/10.1021/i100006a005}
}

@article{Brooks1992SMA,
  author  = {Brooks, Clayton A. and Cramer, Steven M.},
  title   = {Steric mass-action ion exchange: Displacement profiles and induced salt gradients},
  journal = {AIChE J.},
  year    = {1992},
  volume  = {38},
  number  = {12},
  pages   = {1969--1978},
  url     = {https://doi.org/10.1002/aic.690381212}
}

@article{MattssonNordstrom2004,
  author  = {Mattsson, Ken and Nordstr{\"o}m, Jan},
  title   = {Summation by parts operators for finite difference approximations of second derivatives},
  journal = {J. Comput. Phys.},
  year    = {2004},
  volume  = {199},
  number  = {2},
  pages   = {503--540},
  url     = {https://doi.org/10.1016/j.jcp.2004.03.001}
}

@article{MattssonAlmquist2018BoundaryOptimized,
  author  = {Mattsson, Ken and Almquist, Martin and van der Weide, Edwin},
  title   = {Boundary optimized diagonal-norm {SBP} operators},
  journal = {J. Comput. Phys.},
  year    = {2018},
  volume  = {374},
  pages   = {1261--1266},
  url     = {https://doi.org/10.1016/j.jcp.2018.06.010}
}

@article{SvardNordstrom2006Accuracy,
  author  = {Sv{\"a}rd, Magnus and Nordstr{\"o}m, Jan},
  title   = {On the order of accuracy for difference approximations of initial-boundary value problems},
  journal = {J. Comput. Phys.},
  year    = {2006},
  volume  = {218},
  number  = {1},
  pages   = {333--352},
  url     = {https://doi.org/10.1016/j.jcp.2006.02.014}
}

@article{Svard2014SBPReview,
  author  = {Sv{\"a}rd, Magnus and Nordstr{\"o}m, Jan},
  title   = {Review of summation-by-parts schemes for initial--boundary-value problems},
  journal = {J. Comput. Phys.},
  year    = {2014},
  volume  = {268},
  pages   = {17--38},
  url     = {https://doi.org/10.1016/j.jcp.2014.02.031}
}

@book{HesthavenWarburton2008NodalDG,
  author    = {Hesthaven, Jan S. and Warburton, Tim},
  title     = {Nodal Discontinuous {G}alerkin Methods: Algorithms, Analysis, and Applications},
  publisher = {Springer},
  address   = {New York, NY},
  year      = {2008},
  series    = {Texts in Applied Mathematics},
  url       = {https://doi.org/10.1007/978-0-387-72067-8}
}

@article{Ma2021ADSensitivityComparison,
  author  = {Ma, Yingbo and Dixit, Vaibhav and Innes, Mike and Guo, Xingjian and Rackauckas, Christopher},
  title   = {A Comparison of Automatic Differentiation and Continuous Sensitivity Analysis for Derivatives of Differential Equation Solutions},
  journal = {arXiv preprint arXiv:1812.01892},
  year    = {2021},
  url     = {https://doi.org/10.48550/arXiv.1812.01892}
}

@article{Leweke2018CADET,
  author  = {Leweke, Samuel and von Lieres, Eric},
  title   = {Chromatography Analysis and Design Toolkit ({CADET})},
  journal = {Comput. Chem. Eng.},
  year    = {2018},
  volume  = {113},
  pages   = {274--294},
  url     = {https://doi.org/10.1016/j.compchemeng.2018.02.025}
}

@article{Puttmann2013GRMSensitivities,
  author  = {P{\"u}ttmann, Andreas and Schnittert, Sebastian and Naumann, Uwe and von Lieres, Eric},
  title   = {Fast and accurate parameter sensitivities for the general rate model of column liquid chromatography},
  journal = {Comput. Chem. Eng.},
  year    = {2013},
  volume  = {56},
  pages   = {46--57},
  url     = {https://doi.org/10.1016/j.compchemeng.2013.04.021}
}

@article{Puttmann2016AlgorithmicDifferentiation,
  author  = {P{\"u}ttmann, Andreas and Schnittert, Sebastian and Leweke, Samuel and von Lieres, Eric},
  title   = {Utilizing algorithmic differentiation to efficiently compute chromatograms and parameter sensitivities},
  journal = {Chem. Eng. Sci.},
  year    = {2016},
  volume  = {139},
  pages   = {152--162},
  url     = {https://doi.org/10.1016/j.ces.2015.08.050}
}

@article{Hahn2014AdjointChromatography,
  author  = {Hahn, Tobias and Sommer, Anja and Osberghaus, Anna and Heuveline, Vincent and Hubbuch, J{\"u}rgen},
  title   = {Adjoint-based estimation and optimization for column liquid chromatography models},
  journal = {Comput. Chem. Eng.},
  year    = {2014},
  volume  = {64},
  pages   = {41--54},
  url     = {https://doi.org/10.1016/j.compchemeng.2014.01.013}
}

@article{Meyer2023IndustrialIEXDG,
  author  = {Meyer, Kristian and Ibsen, Mikkel S. and Vetter-Joss, Lisa and Hansen, Ernst B. and Abildskov, Jens},
  title   = {Industrial ion-exchange chromatography development using discontinuous {G}alerkin methods coupled with forward sensitivity analysis},
  journal = {J. Chromatogr. A},
  year    = {2023},
  volume  = {1689},
  pages   = {463741},
  url     = {https://doi.org/10.1016/j.chroma.2022.463741}
}

@article{Meyer2018NodalDGChromatography,
  author  = {Meyer, Kristian and Huusom, Jakob K. and Abildskov, Jens},
  title   = {High-order approximation of chromatographic models using a nodal discontinuous {G}alerkin approach},
  journal = {Comput. Chem. Eng.},
  year    = {2018},
  volume  = {109},
  pages   = {68--76},
  url     = {https://doi.org/10.1016/j.compchemeng.2017.10.023}
}

@article{Meyer2020ChromaTech,
  author  = {Meyer, Kristian and Leweke, Samuel and von Lieres, Eric and Huusom, Jakob K. and Abildskov, Jens},
  title   = {{ChromaTech}: A discontinuous {G}alerkin spectral element simulator for preparative liquid chromatography},
  journal = {Comput. Chem. Eng.},
  year    = {2020},
  volume  = {141},
  pages   = {107012},
  url     = {https://doi.org/10.1016/j.compchemeng.2020.107012}
}

@article{Breuer2023DGSEMCADET,
  author  = {Breuer, Jan M. and Leweke, Samuel and Schm{\"o}lder, Johannes and Gassner, Gregor and von Lieres, Eric},
  title   = {Spatial discontinuous {G}alerkin spectral element method for a family of chromatography models in {CADET}},
  journal = {Comput. Chem. Eng.},
  year    = {2023},
  volume  = {177},
  pages   = {108340},
  url     = {https://doi.org/10.1016/j.compchemeng.2023.108340}
}

@article{Bezanson2017Julia,
  author  = {Bezanson, Jeff and Edelman, Alan and Karpinski, Stefan and Shah, Viral B.},
  title   = {Julia: A Fresh Approach to Numerical Computing},
  journal = {SIAM Rev.},
  year    = {2017},
  volume  = {59},
  number  = {1},
  pages   = {65--98},
  url     = {https://doi.org/10.1137/141000671}
}

@inproceedings{Lattner2004LLVM,
  author    = {Lattner, Chris and Adve, Vikram},
  title     = {{LLVM}: A Compilation Framework for Lifelong Program Analysis \& Transformation},
  booktitle = {International Symposium on Code Generation and Optimization (CGO '04)},
  year      = {2004},
  pages     = {75--86},
  publisher = {IEEE},
  url       = {https://doi.org/10.1109/CGO.2004.1281665}
}

@article{Kobl2024OligoIEX,
  author  = {Kobl, Kilian and Nicoud, Lucr{\`e}ce and Nicoud, Edouard and Watson, Anna and Andrews, John and Wilkinson, Edward A. and Shahid, Muhid and McKay, Christopher and Andrews, Benjamin I. and Omer, Batool A. and Narducci, Olga and Masson, Edward and Davies, Suzanne H. and Vandermeersch, Tobias},
  title   = {Oligonucleotide Purification by Ion Exchange Chromatography: A Step-by-Step Guide to Process Understanding, Modeling, and Simulation},
  journal = {Org. Process Res. Dev.},
  year    = {2024},
  volume  = {28},
  number  = {7},
  pages   = {2569--2589},
  url     = {https://doi.org/10.1021/acs.oprd.4c00013}
}

@book{Nicoud2015ChromProcesses,
  author    = {Nicoud, Roger-Marc},
  title     = {Chromatographic Processes: Modeling, Simulation, and Design},
  publisher = {Cambridge University Press},
  address   = {Cambridge, UK},
  year      = {2015},
  series    = {Cambridge Series in Chemical Engineering},
  url       = {https://doi.org/10.1017/CBO9781139998284}
}

@article{BassiRebay1997,
  author  = {Bassi, F. and Rebay, S.},
  title   = {A high-order accurate discontinuous finite element method for the numerical solution of the compressible {N}avier--{S}tokes equations},
  journal = {J. Comput. Phys.},
  year    = {1997},
  volume  = {131},
  number  = {2},
  pages   = {267--279},
  url     = {https://doi.org/10.1006/jcph.1996.5572}
}

@article{SCHENK2004475,
title = {Solving unsymmetric sparse systems of linear equations with PARDISO},
journal = {Future Gener. Comput. Syst.},
volume = {20},
number = {3},
pages = {475--487},
year = {2004},
url = {https://doi.org/10.1016/j.future.2003.07.011},
author = {Schenk, Olaf and Gärtner, Klaus},
}

@article{Schenk2006,
author = {Schenk, Olaf and Gärtner, Klaus},
journal = {ETNA},
pages = {158--179},
publisher = {Kent State University, Department of Mathematics and Computer Science},
title = {On fast factorization pivoting methods for sparse symmetric indefinite systems.},
url = {http://eudml.org/doc/127439},
volume = {23},
year = {2006},
}

@article{Davis2004UMFPACK,
author = {Davis, Timothy A.},
title = {Algorithm 832: UMFPACK V4.3---an unsymmetric-pattern multifrontal method},
year = {2004},
volume = {30},
number = {2},
url = {https://doi.org/10.1145/992200.992206},
journal = {ACM Trans. Math. Softw.},
pages = {196--199},
numpages = {4},
}

@article{Revels2016ForwardDiff,
  author  = {Revels, Jarrett and Lubin, Miles and Papamarkou, Theodore},
  title   = {Forward-Mode Automatic Differentiation in {J}ulia},
  journal = {arXiv preprint arXiv:1607.07892},
  year    = {2016},
  url     = {https://arxiv.org/abs/1607.07892}
}

@article{Rackauckas2017DiffEq,
  author  = {Rackauckas, Christopher and Nie, Qing},
  title   = {{DifferentialEquations.jl} -- {A} Performant and Feature-Rich Ecosystem for Solving Differential Equations in {J}ulia},
  journal = {J. Open Res. Softw.},
  year    = {2017},
  volume  = {5},
  number  = {1},
  pages   = {15},
  url     = {https://doi.org/10.5334/jors.151}
}

@article{Steinebach2023Rodas5P,
  author  = {Steinebach, Gerd},
  title   = {Construction of {Rosenbrock--Wanner} method {Rodas5P} and numerical benchmarks within the {J}ulia {D}ifferential {E}quations package},
  journal = {BIT Numer. Math.},
  year    = {2023},
  volume  = {63},
  number  = {2},
  pages   = {27},
  url     = {https://doi.org/10.1007/s10543-023-00967-x}
}

@article{Shampine1997MATLABODESuite,
  author  = {Shampine, Lawrence F. and Reichelt, Mark W.},
  title   = {The {MATLAB} {ODE} Suite},
  journal = {SIAM J. Sci. Comput.},
  year    = {1997},
  volume  = {18},
  number  = {1},
  pages   = {1--22},
  url     = {https://doi.org/10.1137/S1064827594276424}
}

@article{KennedyCarpenter2019DIRKStiff,
  author  = {Kennedy, Christopher A. and Carpenter, Mark H.},
  title   = {Diagonally implicit {Runge--Kutta} methods for stiff {ODEs}},
  journal = {Appl. Numer. Math.},
  year    = {2019},
  volume  = {146},
  pages   = {221--244},
  url     = {https://doi.org/10.1016/j.apnum.2019.07.008}
}

@article{Kvaerno2004ESDIRK,
  author  = {Kv{\ae}rn{\o}, Anne},
  title   = {Singly diagonally implicit {Runge--Kutta} methods with an explicit first stage},
  journal = {BIT Numer. Math.},
  year    = {2004},
  volume  = {44},
  number  = {3},
  pages   = {489--502},
  url     = {https://doi.org/10.1023/B:BITN.0000046811.70614.38}
}

@article{Davis2010KLU,
  author  = {Davis, Timothy A. and Natarajan, Ekanathan Palamadai},
  title   = {Algorithm 907: {KLU}, A Direct Sparse Solver for Circuit Simulation Problems},
  journal = {ACM Trans. Math. Softw.},
  year    = {2010},
  volume  = {37},
  number  = {3},
  articleno = {36},
  pages   = {1--17},
  url     = {https://doi.org/10.1145/1824801.1824814}
}

@article{vonLieres2010CADET,
  author  = {von Lieres, Eric and Andersson, Joel},
  title   = {A fast and accurate solver for the general rate model of column liquid chromatography},
  journal = {Comput. Chem. Eng.},
  year    = {2010},
  volume  = {34},
  number  = {8},
  pages   = {1180--1191},
  url     = {https://doi.org/10.1016/j.compchemeng.2010.03.008}
}

@article{Danckwerts1953,
  author  = {Danckwerts, Peter V.},
  title   = {Continuous flow systems. {D}istribution of residence times},
  journal = {Chem. Eng. Sci.},
  year    = {1953},
  volume  = {2},
  number  = {1},
  pages   = {1--13},
  url     = {https://doi.org/10.1016/0009-2509(53)80001-1}
}

@book{Guiochon2006Fundamentals,
  author    = {Guiochon, Georges and Felinger, Attila and Shirazi, Dean G. and Katti, Anita M.},
  title     = {Fundamentals of Preparative and Nonlinear Chromatography},
  edition   = {Second},
  publisher = {Elsevier},
  address   = {San Diego, CA},
  year      = {2006},
}

@article{Hindmarsh2005SUNDIALS,
  author  = {Hindmarsh, Alan C. and Brown, Peter N. and Grant, Keith E. and Lee, Steven L. and Serban, Radu and Shumaker, Dan E. and Woodward, Carol S.},
  title   = {{SUNDIALS}: Suite of Nonlinear and Differential/Algebraic Equation Solvers},
  journal = {ACM Trans. Math. Softw.},
  year    = {2005},
  volume  = {31},
  number  = {3},
  pages   = {363--396},
  url     = {https://doi.org/10.1145/1089014.1089020}
}

@book{GeorgeLiu1981Sparspak,
  author    = {George, Alan and Liu, Joseph W.},
  title     = {Computer Solution of Large Sparse Positive Definite Systems},
  publisher = {Prentice-Hall},
  address   = {Englewood Cliffs, NJ},
  year      = {1981},
}

@article{BenchmarkTools2016,
  author = {Jiahao Chen and Jarrett Revels},
  title = {Robust benchmarking in noisy environments},
  journal = {arXiv preprint arXiv:1608.04295},
  year = {2016},
  url = {https://doi.org/10.48550/arXiv.1608.04295}
}

@article{Montoison2025SparseMatrixColorings,
  author  = {Montoison, Alexis and Dalle, Guillaume and Gebremedhin, Assefaw},
  title   = {Revisiting Sparse Matrix Coloring and Bicoloring},
  journal = {arXiv preprint arXiv:2505.07308},
  year    = {2025},
  url     = {https://arxiv.org/abs/2505.07308}
}

@article{Dongarra1990Level3BLAS,
  author  = {Dongarra, Jack J. and Du Croz, Jeremy and Hammarling, Sven and Duff, Iain S.},
  title   = {A Set of Level 3 Basic Linear Algebra Subprograms},
  journal = {ACM Trans. Math. Softw.},
  year    = {1990},
  volume  = {16},
  number  = {1},
  pages   = {1--17},
  url     = {https://doi.org/10.1145/77626.79170}
}

@article{NORDSTROM2003453,
title = {Finite volume methods, unstructured meshes and strict stability for hyperbolic problems},
journal = {Applied Numerical Mathematics},
volume = {45},
number = {4},
pages = {453--473},
year = {2003},
url = {https://doi.org/10.1016/S0168-9274(02)00239-8},
author = {Nordstr{\"o}m, Jan and Forsberg, Karl and Adamsson, Carl and Eliasson, Peter},
}

@article{Xu2018SpectralHP,
  author  = {Xu, Hui and Cantwell, Chris D. and Monteserin, Carlos and Eskilsson, Claes and Engsig-Karup, Allan P. and Sherwin, Spencer J.},
  title   = {Spectral/hp element methods: Recent developments, applications, and perspectives},
  journal = {J. Hydrodyn.},
  year    = {2018},
  volume  = {30},
  number  = {1},
  pages   = {1--22},
  url     = {https://doi.org/10.1007/s42241-018-0001-1}
}

@article{Moura2015LinearDispersionDiffusion,
  author  = {Moura, Rodrigo C. and Sherwin, Spencer J. and Peir{\'o}, Joaquim},
  title   = {Linear dispersion--diffusion analysis and its application to under-resolved turbulence simulations using discontinuous {G}alerkin spectral/$hp$ methods},
  journal = {J. Comput. Phys.},
  year    = {2015},
  volume  = {298},
  pages   = {695--710},
  url     = {https://doi.org/10.1016/j.jcp.2015.06.020}
}

@article{AlbinKlarmann2016DiagonalNormSBP,
  author  = {Albin, Nathan and Klarmann, Joshua},
  title   = {An Algorithmic Exploration of the Existence of High-Order Summation by Parts Operators with Diagonal Norm},
  journal = {J. Sci. Comput.},
  year    = {2016},
  volume  = {69},
  number  = {2},
  pages   = {633--650},
  url     = {https://doi.org/10.1007/s10915-016-0211-z}
}

\end{document}